\documentclass[11pt]{article}
\usepackage{amsmath}
\usepackage{amssymb}
\usepackage{graphics}
\usepackage{color}
\def\disp{\displaystyle}

\def\tto{\;{\lower 1pt \hbox{$\rightarrow$}}\kern -10pt
\hbox{\raise 2pt \hbox{$\rightarrow$}}\;}

\def\Tilde{\widetilde}
\def\Bar{\overline}

\def\h{\hfill\Box}
\def\R{I\!\!R}
\def\N{I\!\!N}

\def\gph{\mbox{\rm gph}\,}

\def\Min{\mbox{\rm Min}\,}
\def\dom{\mbox{\rm dom}\,}

\def\cl*co{\mbox{\rm cl}^*\mbox{\rm co}\,}
\def\cl{\mbox{\rm cl}\,}

\def\cl{\mbox{\rm cl}\,}

\def\h{\hfill\triangle}

\def\O{\Omega}
\def\Th{\Theta}

\def\emp{\emptyset}

\def\lm{\lambda}

\def\Th{\Theta}

\def\N{I\!\!N}

\def\vzero{{\bf 0}}

\begin{document}
\begin{center}
\vspace*{0.3in} {\bf VARIATIONAL PRINCIPLES IN MODELS OF BEHAVIORAL SCIENCES}\\[3ex]
T. Q. BAO\footnote{Department of Mathematics $\&$ Computer Science, Northern Michigan University, Marquette, Michigan 49855, USA (btruong@nmu.edu).}, B. S. MORDUKHOVICH\footnote{Department of Mathematics, Wayne State University, Detroit, Michigan, USA (boris@math.wayne.edu). Research of this author was partially supported by the USA National Science Foundation under grant DMS-1007132.} and A. SOUBEYRAN\footnote{Aix-Marseille University (Aix-Marseille School of Economics), CNRS \& EHESS, Marseille 13002, France (antoine.soubeyran@gmail.com).}
\end{center}
\noindent{\small{\bf Abstract.} This paper develops some mathematical models arising in behavioral sciences, particularly in psychology, which are formalized via general preferences with variable ordering structures. Our considerations are based on the recent ``variational rationality approach" that unifies numerous theories in different branches of behavioral sciences by using, in particular, worthwhile change and stay dynamics and variational traps. In the mathematical framework of this approach, we derive a new variational principle, which can be viewed as an extension of the Ekeland variational principle to the case of set-valued mappings on quasimetric spaces with cone-valued ordering variable structures. Such a general setting is proved to be appropriate for broad applications to the functioning of goal systems in psychology, which are developed in the paper. In this way we give a certain answer to the following striking question: in the world, where all things change (preferences, motivations, resistances, etc.), where goal systems drive a lot of entwined course pursuits between means and ends---what can stay fixed for a while? The obtained mathematical results and new insights open the door to developing powerful models of adaptive behavior, which strongly depart from pure static general equilibrium models of the Walrasian type that are typical in economics.\vspace*{.1in}

\vspace*{.1in}

}
\newtheorem{Theorem}{Theorem}[section]
\newtheorem{Proposition}[Theorem]{Proposition}
\newtheorem{Remark}[Theorem]{Remark}
\newtheorem{Lemma}[Theorem]{Lemma}
\newtheorem{Corollary}[Theorem]{Corollary}
\newtheorem{Definition}[Theorem]{Definition}
\newtheorem{Example}[Theorem]{Example}
\newenvironment{Proof}{\noindent {\it Proof.\, }\hspace{1pt}}
{\hspace{1pt}\hfill $\triangle$\medskip}
\normalsize

\section{Introduction and Motivations}

In this introductory section, we first describe the major features of some stability/stay and change dynamical models in behavioral sciences and the essence of the ``variational rationality" approach to them. Then we show the need for new mathematical developments concerning variational principles and tools of variational analysis for valuable applications to such models. Finally, we discuss the main goals and contributions of this paper, from both viewpoints of mathematics and applications.

{\bf Stability/Stay and Change Dynamics in Behavioral Sciences.} Recent developments on the modeling in various branches of behavioral sciences (including artificial intelligence, economics, management sciences, decision processes, philosophy, political sciences, psychology, sociology) mainly focus on the
functioning/behavioral dynamics of agents, groups, and organizations. Analyzing these models, two very simple observations come to mind. First, all these disciplines, except static models in microeconomics via the classical Walrasian general equilibrium approach \cite{w74}, advocate that human behaviors are driven by {\em adaptive processes}. Second, the vast majority of models in these areas (called sometimes ``theories of stability/stay and change"), advocate that we live in a world where at the same time many things ``stay" (e.g., habits and routines, equilibria, traps, etc.) while many other things ``change" (e.g., creations, destructions, learning, innovation, attitudes as well as beliefs formation and revision, self regulation, goal setting, goal striving and revision, breaking and forming habits, etc.). As stated by Bridges \cite{b09}, ``we are always stuck in the middle between a current status quo position and future ends."

For the reader convenience, we present in Appendix~3 below a brief survey on stability and changes theories. This may help convincing the reader that dynamical  models inherent in behavior sciences are essentially different from more traditional static equilibrium models of the Walrasian type, and thus they require developing appropriate tools of analysis.

{\bf Variational Rationality in Behavioral Sciences.} Since in behavioral sciences all things change (as is often said: ``the only thing that does not change is change itself"), the main question in models of stability and change dynamics is: why, where, how, and when behavioral processes stop or start to change and how transitions work. To describe these issues, Soubeyran \cite{s09,s10} introduced two main variational concepts: {\em worthwhile changes} and {\em variational
traps} as the end points of a succession of worthwhile single changes. The notion of variational traps includes both aspiration points and equilibria and, roughly speaking, reflects the following. Starting from somewhere and not being precisely in a trap, agents want and try approaching such traps in some feasible and acceptable ways (in the case of aspiration points), while being there, prefer to stay than to move away (in the case of equilibria). The notion is crucial in the {\em variational rationality approach} to modelize human behavior suggested in \cite{s09,s10}. This approach helps us to answer the aforementioned  main question as well as to unify and modelize various theories of stability/stay and change. It shows how to model the course of human activities as a succession of worthwhile changes and stays, i.e., a succession of actions balancing at each step between the following:\vspace*{0.05in}

{\bf (i)} {\em Motivation to change} involving the utility/pleasure of advantages to change, where these advantages represent the difference between the future payoff generated by a new action and the future payoff generated by the repetition of the past action.

{\bf (ii)} {\em Resistance to change} involving the desutility/pain of inconveniences to change, where these inconveniences are the difference between costs to be able to change and costs to be able to stay.\vspace*{0.05in}

All these concepts, including those of actions, states, transitions, means (resources and capabilities), ends (performances, payoffs, intentions, goals, desires, preferences and values), judgments, attitudes and beliefs, require lengthy and quite intricate discussions to be fully justified in each different discipline, which have its particular points of view.

At each step we say that changes are {\em worthwhile} if the motivation to change is larger than a chosen fraction of resistance to change. This fraction
represents an adaptive satisficing-sacrificing ratio, which helps us to choose at each step the current level of satisfaction or accepted sacrifice. As argued by Simon \cite{s55}, being ``bounded rational," the agent is not supposed to optimize during the transition even if he/she can or cannot reach the optimum at the
end. The primary aim of the variational rationality approach is to examine the more or less worthwhile to change transition, which can lead to the desired end/goal points via a succession of worthwhile changes and stays. The major questions are as follows:

{\bf(a)} When do such processes make small steps, have finite length, converge?

{\bf(b)} What is the speed of convergence?

{\bf (c)} Do such processes converge in finite time?

{\bf(d)} For which initial points do they converge?

{\bf (e)} Are such processes efficient, i.e, what are the characteristics of end points, which may be critical points, optima, equilibria, Pareto solutions, fixed points, traps, and others?\vspace*{0.05in}

These questions become {\em mathematical} provided that adequate mathematical models within variational rationality approach are created and suitable tools of mathematical analysis are selected. As advocated in the aforementioned papers by Soubeyran, {\em variational analysis}, a relatively new mathematical discipline based on {\em variational principles}, potentially contains an appropriate and powerful machinery to strongly progress in these directions.

{\bf Variational Analysis.} Modern variational analysis has been well recognized as a rapidly developed area of applied mathematics, which is mainly based on {\em variational principles}. It is much related to optimization in a broad sense (being an outgrowth of the classical calculus of variations, optimal control, and mathematical programming), also applying variational principles and optimization techniques to a wide spectrum of problems that may not be of any variational/optimization nature. The reader can find more details on mathematical theories of variational analysis and its many applications in the now classical monograph by Rockafellar and Wets \cite{rw98} as well as in more recent texts by Attouch, Buttazzo and Michaille \cite{abm05}, Borwein and Zhu \cite{bz05}, and the two-volume book by Mordukhovich \cite{m06} with the numerous references therein.

While there are powerful applications of variational analysis to important models in engineering, physics, mechanics, economics\footnote{We particularly refer the reader to the book \cite{m06} and the more recent paper \cite{bm10a} with the vast bibliographies therein for applications of modern techniques of variational analysis and set-valued optimization to models of welfare economics, which are typical in microeconomics modeling.}, etc., not much has been done on applications of variational analysis to psychology and related areas of behavioral science involving human behavior. Within the variational rationality approach, some mathematical results and applications have been recently obtained in the papers \cite{as11,coss13,fls12,ls12,mos11}. However, much more is needed to be done in this direction to capture the {\em dynamical nature} of human behavior reflected in the variational rationality approach. Among the most important dynamical issues, which should be adequately modelized and resolved via appropriate tools of variational analysis, we mention the following settings:\vspace*{0.05in}

{\bf (i)} {\em Periods of the required change} including:

$\bullet$ {\em course of motivation} (e.g., variable preferences, aspirations, hopes, moving goals, goal setting);

$\bullet$ {\em dynamics of resistance to change} (e.g., successive obstacles to overcome, goal striving), which require new concepts of distances, dissimilarity, and
spaces of paths because actions can be defined as succession of operations;

$\bullet$ {\em dynamics of adaptation} concerning mainly self-regulation problems such as feedbacks, goal revision, goal pursuit, etc.\vspace*{0.05in}

{\bf (ii)} {\em Periods, where nothing is required to change}, namely: temporary or permanent ends as optima, stationary, equilibrium points, fixed points, traps, habits, routines, social norms, etc.\vspace*{0.05in}

Having these ``dynamical" issues in mind, we need to revisit available principles and techniques of variational analysis and to develop new mathematical methods and results, which could be applied to solve adaptive dynamic problems arising the aforementioned goals systems of behavioral sciences. Then variational rationality and variational analysis can gain to co-evolve. Variational analysis aims to provide the main tools for the study of variational rationality, which in turn offers a variety of valuable applications for variational analysis in behavioral sciences.

{\bf Main Objectives and Contributions of the Paper.} The primary objective of this paper is to study {\em goal systems} in {\em psychology} by using variational rationality approach and developing an adequate {\em dynamic technique} of variational analysis. To achieve this aim, we establish a new and nontrivial extension of the fundamental {\em Ekeland variational principle} (abbr.\ EVP) to a special class of {\em set-valued} mappings on {\em quasimetric} spaces with cone-valued {\em ordering variable structures}, which becomes the key for our applications to psychology.

The EVP, as first formulated by Ekeland \cite{e72} for extended-real-valued lower semicontinuous functions on metric spaces, is one of the most powerful results of variational analysis and its applications. It is worth mentioning that the original proof in the seminal paper by Ekeland \cite{e74} is complicated and not constructive, involving transfinite induction via Zorn's lemma. The much simplified proof of the EVP, presented in \cite{e79} as a personal communication from Michael Crandall, is remarkable for our purposes, since it is given by a {\em dynamical process} that itself (besides the result) contains significant information for applications to behavioral sciences. However, neither the setting and proof of the latter paper nor their subsequent numerous extensions given in the literature fully fit the main objectives of this paper required by applications to goal systems in psychology. To proceed successfully in this direction, we develop the (dynamical) approach to set-valued extensions of the EVP implemented by Bao and Mordukhovich \cite{bm07,bm10} for mappings in metric spaces with constant Pareto-type ordering preferences to the significantly more involved case of variable ordering structures in quasimetric spaces.

Then we establish valuable applications of the obtained mathematical results to the goal systems in psychology using and enriching the framework of variational rationality approach by Soubeyran \cite{s09,s10}. This allows us, in particular, to shed new light on the explanation, via successions of worthwhile actions and variational traps leading to the underlying dynamical relationships between means and ends in psychological goal systems.

{\bf Organization of the Paper.} The rest of the paper is organized as follows. Section~2 is devoted to the qualitative description and mathematical modeling of the major goal system in psychology from the viewpoint of variational rationality. Besides these issues, we justify here the importance of an appropriate extension of the EVP and the purposes we intend to meet in this way.

Section~3 is pure mathematical containing the formulation and detailed proof of the main mathematical result of this paper, which is the variational principle discussed above. We also present here an important consequence of this result used in what follows.

Section~4 is devoted to the major applications of the developed mathematical theory to the psychological goal system under consideration. Here we present psychological interpretations of the obtained mathematical results and show that they lead us to rather striking psychological conclusions largely discussed in this section with adding more mathematical details.

After concluding remarks in Section~5 and the reference list for the main part of the paper, we presents four appendixes for the reader convenience; each of them has its own references list. Appendix~1 concerns practical means-ends rationality in behavioral sciences. Appendix~2 presents a brief survey of the literature on goal systems in psychology. In Appendix~3  we discuss major references on stability and change dynamics. Finally, Appendix~4 contains additional discussions on the preference change dynamics in behavioral sciences, mainly in psychology.

\section{Goal Systems in Psychology}
\subsection{Formalization of Goal Systems via Means-End Chain}
In what follows, we define a {\em goal system} consisting of {\em four ingredients}; see Appendix~1 and Appendix~3 for more details and discussions.

{\bf (i)} {\em Means} formalized via elements $x\in X$ belonging to the space of means $X$.

{\bf (ii)} {\em Ways} formalized via elements $\omega\in\Omega(x)\subset\overline{\Omega}$ that depend on the given means $x\in X$, where $\Omega(x)$ is a subset of {\em feasible ways} belonging to some space $\overline{\Omega}$.

{\bf (iii)} {\em ``Means-ways of using these means"} pairs formalized as $\phi=(x,\omega)\in X\times\overline{\Omega}=\overline{\Phi}$. Their collection is denoted by $\Phi:=\left\{\phi=(x,\omega)\in\overline{\Phi}|\;\omega\in\Omega(x)\right\}$.

{\bf (iv)} {\em Ends as vectorial payoffs}. Let $P$ be a space of payoffs. These payoffs can be gains $g\in P$ to be increased (e.g., proximal
goals like performances, revenues, profits, utilities, and pleasures as well as distal goals like wishes, desires, and aspirations). These payoffs can also
be costs, unsatisfied needs, desutility, or pains $f\in P$ to be decreased. For instance, $g\in P$ can be a vector of different gains $g=(g^{1},\ldots,g^{m})\in P=\R^{m}$, or can be a vector $f=(f^{1},f^{2},\ldots,f^{m})\in P=\R^{m}$ of unsatisfied needs. We denote by  $g:(x,\omega)\in X\times \Omega(x)\longmapsto g(x,\omega )\in P$ a vectorial {\em payoff function} and by $f:(x,\omega)\in X\times\Omega(x)\longmapsto f(x,\omega)\in P$ a vectorial {\em cost or unsatisfied need function}.

Taking the above into account, {\em goal systems} can be modelized as {\em set-valued mappings} of the following type. For {\em gains} we have the mapping $G(\cdot):x\in X\longmapsto G(x)=\left\{g(x,\omega)|\;\omega\in\Omega(x)\right\}\subset P$ whose values are subsets of payoffs the agent can get given a vector of means $x\in X$. Similarly, for {\em unsatisfied needs} we have the mapping $F(\cdot):x\in X\longmapsto F(x)=\left\{f(x,\omega)|\;\omega\in\Omega(x)\right\}\subset P$ whose values are subsets of unsatisfied needs.\vspace*{0.05in}

The simplest example we can imagine for a goal system is the least interconnected one, where the unique interconnection between goals comes from the {\em resources constraint} {\bf(3)} described below. To proceed, consider the following data involving $j=1,\ldots,m$ activities:

{\bf(1)} $x\in X=\R^{d}$ is a vector of means to be chosen first.

{\bf(2)} $\omega=(\omega^{1},\ldots,\omega^{j},\ldots,\omega^{m})$ is an allocation of the given means $x$, where $\omega^{j}\in\R^{d}$, $i=1,\ldots,m$, will be chosen later.

{\bf(3)} $\omega^{1}+\ldots+\omega^{j}+\ldots+\omega^{m}=x$ is a resource constraint. It defines the way in which the agent allocates the given means $x$ to each activity, namely: the different allocations of means, which can be identified to ways of using means, $\omega^{j}\in X$, aim to reach the goal $g^{j}$ in the activity $j$. It tells us that this allocation is feasible (without slack). This resource constraint can be written in the form
$$
\omega^{1}+\ldots+\omega^{j}+\ldots+\omega^{m}=x\;\Longleftrightarrow\;\omega\in\Omega(x).
$$

{\bf(4)} $g=(g^{1},\ldots,g^{j},\ldots,g^{m})\in P=\R^{m}$ is a vector of goals.

{\bf(5)} $g^{j}=g^{j}(x^{j},\omega^{j})\in\R$ as $j=1,\ldots,m$ represents, relative to the activity $j$, the goal level function $g^{j}(\cdot,\cdot):(x^{j},\omega^{j})\in X\times\overline{\Omega}\longmapsto g^{j}=g^{j}(x^{j},\omega^{j})\in\R$. It tells us that the means $\omega^{j}\in X$ help to reach the goal level $g^{j}=g^{j}(x^{j},\omega ^{j})$. Then $G(x)=\left\{g(x,\omega ),\;\omega\in\Omega(x)\right\}$ defines a goal system as the set-valued ``gain function" $G(\cdot):x\in X\longmapsto G(x)\subset P$. Similarly, $F(x)=\left\{f(x,\omega),\;\omega\in\Omega(x)\right\}$ defines a goal system as the set-valued ``costs or unsatisfied needs function" $F(\cdot):x\in X\longmapsto F(x)\subset P$, where $f=(f^{1},\ldots,f^{j},\ldots,f^{m})\in P=\R^{m}$ and $f^{j}=f^{j}(x^{j},\omega^{j})\in\R$, $j=1,\ldots,m$.

\subsection{Variational Rationality Model of Human Behavior}
{\bf Simplest Adaptive Variational Rationality Model.} The core of the variational rationality approach \cite{s09,s10} can be summarized by the following basic adaptive prototype, which allows a lot of
variants and extensions.\\[1ex]
{\bf (A) Adaptive processes of worthwhile changes and stays.} Agent's behavior is defined as a succession $\left\{x_{0},\ldots,x_{n},\ldots\right\}\subset X$ of actions entwining possible stays $x_{n}\in X\curvearrowright x_{n+1}\in X$, $x_{n+1}=x_{n}$ and possible changes $x_{n}\in X\curvearrowright x_{n+1}\in X$, $x_{n+1}\ne x_{n}$. This behavior is said to be {\em variational rational} if, at each step $n+1$, the agent chooses to change or to stay, depending on what he accepts to consider as the worthwhile change. Then the agent follows a succession of worthwhile stays and changes $x_{n+1}\in W_{\xi_{n+1}}(x_{n})$, $\xi_{n+1}\in\Upsilon$ as $n\in\N$. Let us be more precise.

At step $n$, the agent performs the action $x_{n}$, given the degree of acceptability $\xi_{n}\in\Upsilon$ (to be defined later) he/she has chosen
before. At step $n+1$, given the past action $x_{n}$ done right before and the previously given degree of acceptability $\xi _{n}\in\Upsilon $, the agent
adapts his/her behavior in the following way. He/she chooses a new degree of acceptability $\xi _{n+1}\in\Upsilon$ (which can be the same as before)
of a next worthwhile change $x_{n+1}\in W_{\xi_{n+1}}(x_{n})$. This degree of acceptability (satisficing with some tolerable sacrifices) represents how much worthwhile the agent considers that a change must be to accept to change this step, rather than to stay. There are two cases:

{\bf (i)} A {\em temporary worthwhile stay} $x_{n}\curvearrowright x_{n+1}=x_{n}$. It is the case when $W_{\xi_{n+1}}(x_{n})=\left\{ x_{n}\right\}$. Then the agent will choose, in a rational variational way, to stay at $x_{n}=x_{n+1}$ this time. If at the  next steps $n+2,n+3,\ldots$, the agent does not change the degree of acceptability, he/she will choose to stay there forever. This defines a ``worthwhile to stay" trap, which is a {\em permanent worthwhile stay}.

{\bf (ii)} {\em A temporary worthwhile change} $x_{n}\curvearrowright x_{n+1}\ne x_{n}$. It is the case if $W_{\xi _{n+1}}(x_{n})\ne\left\{x_{n}\right\}$ and if the agent can find $x_{n+1} \in W_{\xi_{n+1}}(x_{n})$ with $x_{n+1}\ne x_{n}$. Then the agent will choose to move from $x_{n}$ to $x_{n+1}\in W_{\xi_{n+1}}(x_{n})$, and so on.\\[1ex]
{\bf(B) Transition phase: the definition of a worthwhile to change step}. Consider step $n+1$, and let $x=x_{n}$ be the preceding action. At step $n+1$, the agent will choose the acceptability ratio $\xi^{\prime}=\xi_{n+1}\in\R_{+}$ and a new action $x^{\prime}=x_{n+1}$. Let $M(x,x^{\prime})\in\R$ be the motivation to change from $x$ to $x^{\prime}$, and let $R(x,x^{\prime})\in\R_{+}$ be the resistance to change from $x$ to $x^{\prime}$. Then the agent will consider that, from his/her point of view, it is worthwhile to move from $x$ to $x^{\prime}$ if the agent's motivation to change is bigger than his/her resistance to change up to the acceptability ratio $\xi_{n+1}$, i.e., under the validity of the condition $M(x,x^{\prime})\ge\xi^{\prime}R(x,x^{\prime})$.

{\em Motivation to change} $M(x,x^{\prime})=U\left[ A(x,x^{\prime})\right]$ is defined as the {\em pleasure} or {\em utility} $U\left[A\right]$ of the advantage to change $A(x,x^{\prime})\in\R$ from $x$ to $x^{\prime}$. In the simplest (separable) case, {\em advantages to change} are defined as the difference $A(x,x^{\prime})=g(x^{\prime})-g(x)$ between a payoff to be improved (e.g., performance, revenue, profit) $g(x^{\prime})\in\R$ when the agent performs a new action $x^{\prime}$ and the payoff $g(x)\in\R$ when he/she repeats a past action $x$ supposing that repetition gives the same payoff as before. On the other hand, advantages to change $A(x,x^{\prime})=f(x)-f(x^{\prime})$ can also be the difference between a payoff $f(x)$ to be decreased (e.g., cost, unsatisfied need) when the agent repeats the same old action $x$ and the payoff $f(x^{\prime})$ the agent gets when he/she performs a new action $x^{\prime}$. The {\em pleasure function} $U\left[\cdot\right]:A\in\R\longmapsto U\left[A\right]\in\R$ is strictly increasing with the initial condition $U\left[0\right]=0$.

{\em Resistance to change} $R(x,x^{\prime})=D\left[I(x,x^{\prime})\right]$ is defined as the {\em pain} or {\em disutility} $D\left[I\right]$ of the
inconveniences to change $I(x,x^{\prime})=C(x,x^{\prime})-C(x,x)\in\R_{+}$, which is the difference between the costs to be able to change $C(x,x^{\prime})\in\R_{+}$ from $x$ to $x^{\prime}$ and the costs $C(x,x)\in\R_{+}$ to be able to stay at $x$. In the simplest case, costs to be able to stay are
supposed to be zero, $C(x,x)=0$ for all $x\in X$, and costs to be able to change are defined as the {\em quasidistances} $C(x,x^{\prime})=q(x,x^{\prime})\in\R_{+}$ satisfying: {\bf(a)} $q(x,x^{\prime})\ge 0$, {\bf (b)} $q(x,x^{\prime})=0\Longleftrightarrow x^{\prime }=x$, and {\bf(c)} $q(x,x'')\le q(x,x^{\prime})+q(x^{\prime},x'')$ for all $x,x^{\prime},x''\in X$. The {\em pain function} $D\left[\cdot\right]:I\in\R_{+}\longmapsto D\left[I\right]\in\R_{+}$ is strictly increasing with the initial condition $D\left[0\right]=0$.

Thus, in this simplest case, the worthwhile to change and stay process satisfies the conditions $g(x_{n+1})-g(x_{n})\ge\xi_{n+1}q(x_{n},x_{n+1})$ or $f(x_{n})-f(x_{n+1})\ge\xi_{n+1}q(x_{n},x_{n+1})$ for each $n$. This yields
$$
W_{\xi^{\prime}}(x)=\left\{x^{\prime }\in X\big|\;g(x^{\prime})-g(x)\ge\xi^{\prime}q(x,x^{\prime})\right\}\;\mbox{ or }\;W_{\xi^{\prime}}(x)=\left\{x^{\prime}\in X\big|\;f(x)-f(x^{\prime})\ge\xi^{\prime }q(x,x^{\prime})\right\}.
$$
{\bf(C) End points as traps.} Given the final worthwhile to change rate $\xi_{\ast}>0$, we say that the end point $x_{\ast}\in X$ of the process under consideration is a {\em stationary trap} if $W_{\xi_{\ast}}(x_{\ast})=\left\{x_{\ast}\right\}$. In the simplest case above, $x_{\ast}\in X$ is a stationary trap if $g(x^{\prime})-g(x_{\ast})<\xi_{\ast}q(x_{\ast},x^{\prime})$ or $f(x_{\ast})-f(x^{\prime})<\xi _{\ast}q(x_{\ast},x^{\prime})$ for all $x^{\prime}\ne x_{\ast}$. The definition of a {\em variational trap} requires more. It involves the given initial state $x_0$ and requires that the final stationary state (trap) can be reachable from this initial state via a worthwhile and feasible transition of single worthwhile changes and temporary stays.\\[1ex]
{\bf(D) Variational rationality problems} include the following major components.

Starting from any given $x_{0}\in X$ and depending on the motivation and resistance to change functions, we want to find a {\em path of worthwhile changes} so that:

{\bf(i)} the steps go to {\em zero} and have {\em finite length};

{\bf (ii)} the corresponding iterations {\em converge to a variational trap};

{\bf(iii)} the {\em convergence rate} and {\em stoping criteria} are investigated;

{\bf (iv)} the {\em efficiency} or {\em inefficiency} of such worthwhile to change processes are studied to clarify whether the worthwhile to change process ends at a {\em critical point}, a {\em local or global optimum}, a {\em local or global equilibrium}, an {\em epsilon-equilibrium}, a {\em Pareto solution}, etc.

\subsection{Ekeland's Variational Principle (EVP) in the Simplest Nonadaptive Model of Variational Rationality}
Let us discuss here how the variational rationality approach of \cite{s09,s10} interprets the classical EVP \cite{e74} in the case of the simplest nonadaptive model. First recall the seminal Ekeland's result.\\
\textbf{The classical EVP}. {\em Let $(X,d)$ be a complete metric space, and let $f(\cdot):x\in X\longmapsto f(x)\in\R\cup\left\{\infty\right\}$ be a lower semicontinuous (l.s.c.) function not identically to $\infty$ and bounded from below. Denote $\underline{f}:=\inf\left\{f(x)|\;x\in X\right\}>-\infty$. Then for every
$\varepsilon>0$, $\lambda>0$, and $x_{0}\in X$ with $f(x_{0})<\underline{f}+\varepsilon$ there exists $x_{\ast}\in X$ satisfying the conditions:
\begin{itemize}
\item[\bf{(a)}] $f(x_{0})-f(x_{\ast})\ge(\varepsilon/\lambda)d(x_{0},x_{\ast})$,

\item[\bf{(b)}] $f(x_{\ast})-f(x^{\prime})<(\varepsilon/\lambda)d(x_{\ast},x^{\prime})$ for all $x^{\prime}\ne x_{\ast}$,

\item[\bf{(c)}] $d(x_{0},x_{\ast})\le\lambda$.
\end{itemize}}

By taking $f(x):=-g(x)$, we can immediately reformulate the EVP for the case of {\em maximization} of $g(.):x\in X\longmapsto g(x)\in\R\cup\left\{-\infty\right\}$. Let us next present a {\em variational rationality} interpretation of the latter result by using terminology and notation of Sections~2.1 and 2.2.\\[1ex]
{\bf Variational rationality interpretation of the EVP}. Consider the maximization formulation for a payoff to be improved. Then the EVP variational rationality framework tells us the following. Impose the {\bf assumptions}:

$\bullet$ the worthwhile to change process $x_{n+1}\in W_{\xi _{n+1}}(x_{n})$ is {\em nonadaptive}, which means that the ``satisficing-sacrificing" ratio $\xi_{n+1}$ is constant along the process $\xi_{n+1}\equiv\xi=\varepsilon/\lambda>0$
for all $n\in N$;

$\bullet$ advantages to change are {\em separable}, i.e., $A=A(x,x^{\prime})=g(x^{\prime})-g(x)$, where $g(\cdot):x\in X\longmapsto g(x)\in\R$ is a {\em payoff function} to be improved (in the sense of maximization);

$\bullet$ costs to be able to change $C=C(x,x^{\prime})\in\R_{+}$ represent a {\em distance} $C(x,x^{\prime})=d(x,x^{\prime })$, which implies that costs to be
able to change are {\em symmetric} $C(y,x)=C(x,y)$, costs to be able to stay are {\em zero} $C(x,x)=0$ for all $x\in X$, and costs to be able to change satisfy the {\em triangular inequality};

$\bullet$ pleasure and pain are identified with {\em advantages to change} and {\em inconveniences to change}, respectively, i.e., $U\left[A\right]=A$ for all $A\in\R,$ and $D\left[I\right]=I$ for all $I\ge 0$.\\[1ex]
Then we have the {\bf conclusions:}
\begin{itemize}
\item[\bf(a)] There exists an acceptable {\em one step transition} from any initial position $x_0$ to the end $x_{\ast}\in W_{\xi}(x_{0})$.
This means that it is {\em worthwhile} to move directly from $x_{0}$ to $x_{\ast}$.

\item[\bf(b)] The end is a {\em stable position}, which means that $W_{\xi}(x_{\ast})=\left\{x_{\ast}\right\}$. In other words, it  is
{\em not worthwhile} to move from $x_{\ast}$ to any different action $x^{\prime}\ne x_{\ast}$.

\item[\bf(c)] The end can be reached in a {\em feasible way} $C(x_{0},x_{\ast})\le\lambda$. This means that if the agent cannot spend more
than the $\lambda>0$ amount in terms of costs, then the move from $x_{0}$ to $x_{\ast}$ is feasible in the model.
\end{itemize}

\subsection{Variational Traps and Behavioral Essence of the Ekeland Principle}

As state by Alber and Heward \cite{ah96}, the essence of a trap, given in behavioral terms, is that only ``a relatively simple response is necessary to enter the trap, yet once entered, the trap cannot be resisted in creating general behavior changes." Let us give (among many others) a short list of traps we can find in different disciplines.\\[1ex]
{\bf (A)} {\bf Psychology}. Baer and Wolf \cite{bw70} seem to be the first to use the term of behavioral trap in describing ``how natural contingencies of reinforcement operate to promote and maintain generalized behavior changes." Plous \cite{p93} lists five behavioral traps defined as more or less easy to fall into and more or less difficult to get out: investment, deterioration, ignorance, and collective traps. Behavioral traps have been shown to end reinforcement
processes \cite{s04}. Ego-depletion can generate behavioral traps due to fatigue costs, in the context of self regulation failures \cite{b02,bh96}. Among several cognitive and emotional traps we can list all-or-nothing thinking, labeling, overgeneralization, mental filtering, discounting the positive, jumping to conclusions, magnification, emotional reasoning, should and shouldn't statements, personalizing the blame, etc.\\[1ex]
{\bf (B)} {\bf Economics and decision sciences}.  Making traps in decision represents hidden biases, heuristics, and routines; e.g., anchoring, status quo, sunk costs, confirming evidence, framing, estimation, and forecasting traps; see \cite{hkr98} and the references therein.\\[1ex]
{\bf(C)} {\bf Management sciences}. The importance of success and failure traps within organizations due to the so-called ``myopia of learning" is emphasized in \cite{lm93,m91}.\\[1ex]
{\bf (D)} {\bf Development theory}. To explain the formation of poverty traps, Appadurai \cite{a04} defines {\em aspiration traps}, which describe the inability to aspire of the poor; see also \cite{hm06,r06}.\vspace*{0.05in}

{\em Variational approach} of \cite{s09,s10} shows that, from the viewpoint of behavioral sciences dealing with essentially {\em dynamic models} of human behaviors (contrary to pure static developments in general equilibrium theory of economics), the very {\em essence of the EVP} concerns {\em variational traps}. More precisely, conditions {\bf(a)} and {\bf(c)} of the Ekeland theorem presented above define a variational trap, which is rather easy to reach in an acceptable and feasible way, while is difficult to leave due to condition {\bf(b)}. This corresponds to the intuitive sense of variational traps in behavioral sciences given in \cite{p93}. From this viewpoint, the EVP not only ensures the {\em existence} of variational traps, but also indicates (particularly in its proof) the {\em dynamics} of how to reach a variational trap.

It is worth mentioning that the usage and understanding of the EVP in the variational rationality approach to behavioral sciences is different from those in mathematics. Indeed, in {\em behavioral sciences} (where inertia, frictions, and learning play a major role), natural solutions are variational traps that are reachable in a worthwhile way as {\em maximal elements} of certain {\em dynamic} relationships for {\em worthwhile changes}. In this way, the exact solutions become variational traps, since they include costs to be able to change in their definition. The approximate solution becomes optimum, since they ignore costs to be able to change in their definition.

In {\em mathematics}, the treatment of the EVP is actually {\em opposite}. Variational traps resulting from the EVP are seen as {\em approximate solutions} to the original problem while providing the {\em exact optimum} to another optimization problem, with a small {\em perturbation} term.

\subsection{Variationally Rational Model of Goal Systems}
{\bf Variational Rationality Concepts: Worthwhile to Change Payoffs.} In the context of goal systems, we define the following variational
concepts following \cite{s09,s10}.

{\bf 1. Changes.} We say that $\phi=(x,\omega)\in\Phi\mathbf{\curvearrowright}\phi^{\prime}=(x^{\prime},\omega^{\prime})\in\Phi$ signifies a {\em change} from the old feasible ``means-way of using these means" pair $\phi\in\Phi$ to the new feasible pair $\phi^{\prime}\in\Phi$, where
$$
\Phi:=\left\{\phi=(x,\omega)\in\overline{\Phi}\;\mbox{ such that }\;\omega\in\Omega(x)\right\}
$$
stands for the set of all the feasible pairs.

{\bf 2. Advantages to change.} Consider now a change  from the present feasible means-ends pair $(x,g)$ with $g\in G(x)$ to the next one $(x^{\prime},g^{\prime})$ with $g^{\prime}\in G(x^{\prime})$. Then $A:=A(\phi,\phi^{\prime})=g(\phi^{\prime})-g(\phi)\in P$ is the
{\em advantage} to change from the old feasible pair $\phi\in\Phi$ to the new feasible pair $\phi^{\prime}\in\Phi$.

{\bf 3. Costs to be able to change and costs to be able to stay.} Denote by $C(\cdot,\cdot):(\phi,\phi^{\prime})\in\Phi\longmapsto C\left[\phi,\phi^{\prime}\right]\in P$ the {\em costs to be able to change} from the old feasible pair $\phi\in\Phi$ to the new feasible pair $\phi^{\prime }\in\Phi$. It is worth mentioning here that, in the context of our new version of the EVP for {\em variable ordering structures} developed in Section~3, the costs to be able to change exhibit the following two specific properties:

{\bf(i)} They do {\em not depend} on the ways of using means $C\left[\phi,\phi^{\prime}\right]=C(x,x^{\prime})\in P$. This means that they actually behave as if the ways of using means were free.

{\bf(ii)} They have a {\em directional shape} $C(x,x^{\prime})=q(x,x^{\prime})\xi$, where $\xi\in P$ and $q(x,x^{\prime})\in\R_{+}$ is a scalar
{\em quasidistance}, which represents the total cost to be able to change from the old means $x$ to the new means $x^{\prime}$. In the case where
$P=\R^{J}$, the vector $\xi=(\xi^{1},\ldots,\xi^J)\in P$ with $\xi^j\in\R_{+}$, $j=1,\ldots,J$, and $\left\Vert\xi\right\Vert =1$ represents the {\em internal shares} of this scalar total cost $q(x,x^{\prime})$ among different activities. In general these shares $\xi^{j}=\xi^{j}(x,x^{\prime})>0$ can change along the process.

The detailed justification that the total costs to be able to change can be modelized as a quasidistance $q(x,x^{\prime})\in\R_{+}$ is given in \cite{s09}. To save space, let us just mention that this comes from the definition of the costs to be able to change as the infimum of the costs to be
able to perform a succession of operations of deletions, conservations, and acquisitions. The fact that the costs to be able to stay satisfy $C(x,x)=0$ for all $x\in X$ must also be carefully justified. In the general case, the costs to be able to change modelize inertia, i.e., the resistance to change. There are two extreme cases. {\em Strong resistance to change} is modelized by the costs to be able to change as scalars or cone quasidistances. This is the case of {\em variational principles} of Ekeland's type. On the other hands, {\em weak resistance} to change is modelized by the costs to be able to change via {\em convex increasing functions} of scalar or cone quasidistances. This is the case of {\em proximal
algorithms}; see \cite{bs13,mos11} for more details and discussions.

{\bf 4. Inconveniences to change.} They represent the difference $I(\phi,\phi^{\prime})=C(x,x^{\prime})-C(x,x)$ between the costs to be able to change $C(x,x^{\prime })$ and the costs to be able to stay $C(x,x)$.

{\bf 5. Worthwhile to change payoffs.} Consider the difference between the advantages to change and the costs to be able to change given by
$$
\Delta:=\Delta\left[(x,\omega),(x^{\prime},\omega^{\prime})\right]=\Delta\left[\phi,\phi^{\prime}\right]=A(\phi,\phi^{\prime})-\xi I(\phi,\phi^{\prime})=\big(g(\phi^{\prime})-g(\phi)\big)-\xi q(x,x^{\prime})\in P.
$$
This defines the worthwhile to change payoff for the change $\phi=(x,\omega)\curvearrowright\phi^{\prime }=(x^{\prime},\omega^{\prime})$, where $\phi,\phi^{\prime}\in\Phi$. Then the {\em change $\phi:=(x,\omega)\curvearrowright\phi^{\prime}=(x^{\prime},\omega ^{\prime })$ is worthwhile} if $\Delta\left[\phi,\phi^{\prime}\right]\ge_{K\left[f(\phi)\right]}\mathbf{0}$.

{\bf 6. Pleasure and pain.} To simplify our model of goal systems in this paper, we will not consider the pleasure and pain functions in full generality, i.e., defined as the utilities $U\left[A(\phi,\phi^{\prime})\right]\in{\bf U}$ of the advantages to change and the pains $D\left[I(\phi,\phi^{\prime})\right]\in{\bf D}$ as the disutilities of inconveniences to be able to change. We simply identify the pleasures with the advantages to change $U\left[A\right]=A$
and the pains with the inconveniences to be able to change $D\left[I\right]=I$. However, the variable cones $K\left[f(\phi)\right]$ or $K\left[g(\phi)\right]$ represent these variable pleasures and pains feelings. They define {\em variable preferences} in the payoff space $P$. Then the change $\phi=(x,\omega)\curvearrowright\phi^{\prime }=(x^{\prime},\omega^{\prime})$ is {\em worthwhile} if $\Delta\left[\phi,\phi^{\prime}\right]\ge_{K\left[f(\phi)\right]}\mathbf{0}$. This defines the corresponding {\em variable preference on feasible ``means-way of using these means" pairs} in the following way, respectively:
$$
\phi''\ge_{\phi}\phi^{\prime}\Longleftrightarrow\Delta\left[\phi,\phi''\right]\ge_{K\left[f(\phi)\right]}\Delta\left[\phi,\phi^{\prime}\right],
$$
$$\phi''\ge_{\phi}\phi^{\prime}\Longleftrightarrow\Delta\left[\phi,\phi''\right]\geq_{K\left[g(\phi)\right]}\Delta\left[\phi,\phi^{\prime}\right],
$$
where the reference point is the current feasible pair $\phi=(x,\omega)\in\Phi$.

\section{Variational Principle for Variable Ordering Structures}

The preceding section describes in detail the primary adaptive psychological model of this paper and also discusses the importance of {\em variational analysis} (particularly an appropriate variational principle of the Ekeland type) as the main mathematical tool of our study and applications. In comparison with the original version of the EVP presented above, the following {\em three requirements} are absolutely mandatory for an appropriate extension of the EVP for its possible application to the psychological model under consideration:

{\bf (a)} {\em vectorial} (actually {\em set-valued}) nature of the cost function;

{\bf (b)} {\em quasimetric} (instead of metric) structure of the topological space of arguments;

{\bf (c)} {\em variable preference} structure of ordering on the space of values.\vspace*{0.05in}

By now, a great many of numerous versions and extensions of the EVP are known in the literature; see, e.g., \cite{bm07,bm10,bz05,grtz03,ln11,m06,q12} and the references therein for more recent publications. Some of them address the above issues {\bf (a)} and {\bf (b)} while {\em none} of them, to the best of our knowledge, deals with {\em variable structures} in {\bf (c)}, which is the main issue required for applications to {\em adaptive} psychological models as well as to adaptive models in other branches of behavioral sciences.

Note that problems with variable preferences have drawn some attention in recent publications (e.g., \cite{bm13,bb13,cy02,e11,eh13,ls12}), but not from the viewpoint of variational principles as in this paper.\vspace*{0.05in}

In this section, we derive a general variational principle that addresses all the three issues {\bf(a)}--{\bf(c)} listed above. Furthermore, we consider a general {\em parametric} setting when the mapping in the variational principle depends on a {\em control} parameter, which allows us to take into account the ``ways of using means" providing in this way a kind of {\em feedback} in adaptive psychological models; see Section~4 for more discussions. Our approach and main result extend those (even in nonparametric settings of finite-dimensional spaces) from the papers by Bao and Mordukhovich \cite{bm07,bm10}, which dealt with nonparametric mappings between Banach spaces in the standard (not variable) preference framework. Addressing the new challenges in this paper requires a significant improvement of the previous techniques, which is done below.\vspace*{0.05in}

To describe the class of variable preferences invoking in our main result, take vectors $p_1,p_2\in P$ from some linear topological {\em decision space} $P$, denote $d:=p_1-p_2$, and say that $p_2$ is {\em preferred} by the decision maker to $p_1$ with the {\em domination factor} $d$ for $p_1$. The set of all the  domination factors for $p_1$ together with the zero vector $\vzero\in P$ is denoted by $K[p_1]$. Then the set-valued mapping $K:P\tto P$ is called a {\em variable ordering structure}. We define the {\em ordering relation} induced by the variable ordering structure $K$ by
\begin{eqnarray*}
p_2\;\le_{K[p_1]}p_1\;\mbox{ if and only if }\;p_2\in p_1-K[p_1]
\end{eqnarray*}
and say that $p_\ast\in\Xi$ is {\em Pareto efficient/minimal} to the set $\Xi$ in $P$ with respect to the ordering structure $K$ if there is no other vector $p\in \Xi\setminus\{p_\ast\}$ such that $p\le_{K[p_\ast]}p_\ast$, i.e.,
$$
\big(p_\ast-K[p_{\ast}]\big)\cap\Xi=\{p_\ast\}.
$$
It follows from the definition that $p_\ast\in\Min(\Xi;K[p_{\ast}])$ in the sense of set optimization with the ordering cone $K[p_\ast]$; see, e.g., \cite{grtz03}. This order reduces to the one in set optimization when $K[p]\equiv\Theta$ for some convex ordering cone $\Theta\subset P$.

Fixing a {\em direction} $\xi\in P$ and a {\em threshold/accuracy} $\varepsilon>0$, we say that $p_\ast$ is an {\em approximate $\varepsilon\xi$-minimal point} of $\Xi$ with respect to $K$ if
\begin{eqnarray*}
\big(p_{\ast}-K[p_{\ast}]-\varepsilon\xi\big)\cap\Xi=\emp.
\end{eqnarray*}

Next we recall the definition of quasimetric spaces and the corresponding notions of closedness, compactness, and completeness in such topological spaces.

\begin{Definition} {\bf (quasimetric spaces).} A pair $(X,q)$ with the collection of elements $X$ and the function $q:X\times X\longmapsto\R$ on  $X\times X$ is said to be a {\sc quasimetric space} if the following hold:
\begin{itemize}
\item[\bf(i)] $q(x,x^{\prime})\ge 0$ for all $x,x^\prime\in X$;
\item[\bf(ii)] $q(x,x^{\prime})=0$ if and only if $x^{\prime}=x$ for all $x,x^\prime\in X$;
\item[\bf(iii)] $q(x,x'')\le q(x,x^{\prime})+q(x^{\prime },x'')$ for all $x,x^{\prime},x''\in X$.
\end{itemize}
\end{Definition}

In what follows, we consider only those complete quasimetric spaces, where each Cauchy sequence converge to the {\em unique} limit, which is of course automatic in metric spaces. This is a typical setting for applications of quasimetric spaces in behavioral sciences and other disciplines. The reader can find more discussions on such quasimetric spaces and sufficient conditions for the aforementioned uniqueness property in the recent papers \cite{bs13,ft13} and the references therein.

\begin{Definition}{\bf (left-sequential closedness).} A subset $S\subset X$ is said to be {\sc left-sequentially closed} if for any sequence $\left\{x_{n}\right\}\subset X$ converging to $x_{\ast}\in X$ in the sense that the numerical sequence $\left\{q(x_{n},x_{\ast})\right\}$ converges to zero, the limit $x_{\ast}$ belongs to $S$.
\end{Definition}

\begin{Definition}{\bf (left-sequential completeness).} A sequence $\left\{x_{n}\right\}\subset X$ is said to be {\sc left-sequential Cauchy} if for each $k\in\N$ there exists $N_{k}$ such that
$$
q(x_{n},x_{m})<1/k\;\mbox{ for all }\;m\ge n\ge N_{k}.
$$
A quasimetric space $(X,q)$ is said to be {\sc left-sequentially complete} if each left-sequential Cauchy sequence is convergent.
\end{Definition}

Let $f:T\rightarrow P$ be a mapping from a quasimetric space $(T,q)$ to an ordered vector space $P$ equipped with a variable ordering structure $K: P\rightrightarrows P$, and let $S\subset T$. Then:

$\bullet$ $f$ is (left-sequentially) {\em level-closed with respect to} $K$ if for any $p\in P$ the $p$-level set of $f$ with respect to $K$ defined by
\begin{eqnarray*}
\mbox{\rm lev}_{p}(f,K):=\left\{t\in X\big|\;f(t)\le_{K\left[p\right]}p\right\}=\left\{t\in X\big|\;f(t)\in p-K\left[ p\right]\right\}
\end{eqnarray*}
is left-sequentially closed in $(T,q)$.

$\bullet$ $f$ is {\em quasibounded from below} on $S\subset\dom f:=\{t\in T|\;f(t)\in P\}$ {\em with respect to a cone} $\Theta$, or it is {\em $\Theta$-quasibounded from below} for short, if there is a bounded subset $M\subset P$ such that $f(S)\subset M+\Theta$ for the image set $f(S):=\cup\{f(t)\in P|\;t\in S\}$.

$\bullet$ $t^\ast\in S$ is a {\em Pareto minimizer} (resp.\ {\em approximate $\varepsilon\xi$-minimizer}) of $f$ over $S$ with respect to $K$ if $f(t_\ast)$ is the corresponding Pareto minimal point (resp.\ approximate $\varepsilon\xi$-minimal point) of the image set $f(S)\subset P$.\vspace*{0.05in}

Note that our applications in Section~4, $x\in X$ represents actions, states, or some means; $g\in P$ represents vectors of ends to be increased (e.g., performances, payoffs, revenues, profits, utilities, pleasures), and $f\in P$ represents vectors of ends to be decreased (a vector of costs, unsatisfied needs, disutilities, pains, etc.). Until arriving at applications, in the mathematical framework here we consider for definiteness the ``minimization" setting (instead of the ``maximization" one), which is more appropriate for certain applications. Correspondingly, $f\in P$ as a vector of ends to be {\em decreased} and $K\left[f\right]\subset P$ is the cone of vectorial costs {\em lower} than the given cost vector $f$. We say that the vectorial payoff $f_{2}$ (a vector of payoffs to be decreased) is {\em smaller} than $f_{1}$ with respect to $K$ and write $f_{2}\le_{K\left[f_{1}\right]}f_{1}$ if $f_{2}\in f_{1}-K\left[ f_{1}\right]$.\vspace*{0.05in}

Consider now a set-valued mapping $\Omega:X\tto\Bar\Omega$ from a {\em quasimetric} space $(X,q)$ to a {\em compact} subset $\Bar\Omega\subset Y$ of a {\em Banach} space $Y$. Let $P$ be a {\em linear topological} space (of {\em payoffs}) equipped with some variable {\em cone-ordering} structure $K:P\tto P$ (called {\em variable preference on payoffs}), and let $\emp\ne\Th\subset P$ be a cone. Our standing assumptions are as follows.

Now let us formulate our {\em standing assumptions} on the initial data of the problem under consideration needed for the proof of the new variational principle in Theorem~\ref{EVP-VOS}. Recall that a cone $K\subset P$ is {\em proper} if we have $K\ne\left\{\mathbf{0}\right\}$ and $K\ne P$.

{\bf(H1)} The quasimetric space $(X,q)$ is {\em left-sequentially complete}. Furthermore, the quasimetric $q(x,\cdot)$ is  (left-sequentially) {\em l.s.c.} with respect to the second variable for all $x\in X$.

{\bf(H2)} All the values of $K\colon P\tto P$ are {\em closed, convex}, and {\em pointed} subcones of $P$. Furthermore, the {\em common ordering cone} of $K$, denoted by $\Theta_{K}:=\cap_{f\in P}K\left[f\right]$, also has these properties.

{\bf(H3)} The mapping $K\colon P\tto P$ enjoys the {\em transitivity property} in the sense that
$$
\Big(f_{1}\in f_{0}-K\left[f_{0}\right],\;f_{2}\in f_{1}-K\left[f_{1}\right]\Big)\Longrightarrow\Big(f_{2}\in f_{0}- K\left[f_{0}\right]\Big).
$$

{\bf(H4)} The mapping $\Omega:X\rightrightarrows\overline{\Omega}$ is (left-sequentially) {\em closed-graph}.

{\bf(H5)} The cone $\Th\subset P$ is {\em closed} and {\em convex}.\vspace*{0.05in}

It is easy to check that the relation $f_{1}\le_{K\left[f_{0}\right]}f_{0}$ implies that $K\left[f_{1}\right]\subset K\left[f_{0}\right]$ with the equality
$K\left[f_{1}\right]+K\left[f_{0}\right]=K\left[f_{0}\right]$ under assumption {\bf (H2)}.

\begin{Theorem}{\bf (parametric variational principle for mappings with variable ordering).}\label{EVP-VOS} Let $f:X\times\Bar\Omega\longmapsto P$ be a mapping with $\dom f=\gph\Omega$ in the setting described above. In addition to the standing assumptions {\bf(H1)}--{\bf(H5)}, suppose that:

{\bf(A1)} $f$ is quasibounded from below on $\gph\Omega$ with respect to the cone $\Theta$.

{\bf(A2)} $f$ is (left-sequentially) level-closed with respect to $K$ on $\gph\Omega$.

{\bf(A3)} $f(x,\cdot)$ is continuous for each $x\in\dom\Omega$.\\[1ex]
Then for any $\varepsilon>0$, $\lambda>0$, $(x_{0},\omega_{0})\in\gph\Omega$, and $\xi\in\Theta_{K}\setminus(-\Theta-K\left[f_0\right])$ with $\left\Vert\xi\right\Vert=1$ and $f_0:=f(x_{0},\omega_{0})$ there is a pair $(x_{\ast},\omega_{\ast})\in\gph\Omega$ with
$f_{\ast}:=f(x_{\ast},\omega_{\ast})\in\Min\big(F(x_{\ast});K\left[f_{\ast}\right]\big)$ and $F(x_{\ast}):=\cup\{f(x_{\ast},\omega)|\;\omega\in\Omega(x_{\ast})\}$ satisfying the relationships
\begin{eqnarray}\label{0.1}
f_{\ast}+(\varepsilon/\lambda)q(x_{0},x_{\ast})\xi\le_{K\left[f_{0}\right]}f_{0},
\end{eqnarray}
\begin{eqnarray}\label{0.2}
f+(\varepsilon/\lambda)q(x_{\ast},x)\xi\nleq_{K\left[f_{\ast}\right]}f_{\ast}\;\mbox{ for all }\;(x,\omega)\in\gph\Omega\;\mbox{ with }\;f:= f(x,\omega)\ne f_{\ast}.
\end{eqnarray}
If furthermore $(x_{0},\omega_{0})$ is an approximate $\varepsilon\xi$--minimizer of $f$ over $\gph\Omega$ with respect to $K\left[f_{0}\right]$,
then $x_{\ast}$ can be chosen so that in addition to \eqref{0.1} and \eqref{0.2} we have
\begin{eqnarray}\label{0.3}
q(x_{0},x_{\ast})\le\lambda.
\end{eqnarray}
\end{Theorem}
{\bf Proof.} Without loss of generality we assume that $\varepsilon=\lambda=1$. Indeed, the general case can be easily reduced to this special one by
applying the latter to the mapping $\Tilde f(x,\omega):=\varepsilon^{-1}f(x,\omega)$ and the left-sequentially complete quasimetric space $(X,\Tilde q)$
with $\Tilde q(x,x^\prime):=\lambda^{-1}q(x,x^\prime)$.

Define now a set-valued mapping $W:X\times\Bar\Omega\rightrightarrows X$ by
\begin{eqnarray} \label{0.4}
W(x,\omega):=\left\{x^{\prime}\in X\big|\;\exists\;\omega^{\prime }\in\Omega(x^{\prime})\;\mbox{ with }\;f(x^{\prime},\omega^{\prime})+q(x,x^{\prime})\xi \;\le_{K\,\left[f\right]}f(x,\omega)\right\},
\end{eqnarray}
where $f:=f(x,\omega)$. It is easy to check that for such a pair $(x^{\prime },\omega^{\prime})$ satisfying the inequality in (\ref{0.4})
we get, by denoting $f^{\prime}:=f(x^{\prime},\omega^{\prime})$, that
\begin{eqnarray}\label{0.5}
f(x^{\prime},\omega^{\prime})\;\le_{K\left[f\right]}f(x,\omega)\;\mbox{ and }\;K\left[f^{\prime}\right]\subset K\left[f \right]
\end{eqnarray}
under the imposed assumptions for $K$. Indeed, the inclusion in \eqref{0.5} follows directly from the inequality therein while the latter is valid by
\begin{eqnarray*}
&&f(x^{\prime},\omega ^{\prime})+q(x,x^{\prime})\xi\;\le_{K\left[f\right]}f(x,\omega)\\
&\Longleftrightarrow&f(x^{\prime},\omega^{\prime})+q(x,x^\prime )\xi\in f(x,\omega)-K\left[f\right]\\
&\Longleftrightarrow&f(x^{\prime},\omega^{\prime})\in f(x,\omega)-q(x,x^{\prime})\xi-K\left[f\right]\\
&\Longrightarrow&f(x^{\prime},\omega^{\prime})\in f(x,\omega)-K\left[f\right],
\end{eqnarray*}
where the implication holds due to the choice of $\xi\in\Theta_K\subset K\left[f\right]$ and the convexity of the cone $K\left[f\right]$. In fact $K\left[f\right]+ q(x,x^{\prime})\xi\subset K\left[f\right]+K\left[f\right]=K\left[f\right]$.\vspace*{.05in}

Next we list some important properties of the set-valued mapping $W$ used in what follows.\vspace*{0.05in}

$\bullet$ The sets $W(x,\omega )$ are {\em nonempty} for all $(x,\omega)\in\gph\Omega $ due to $(x,\omega)\in W(x,\omega )$.\vspace*{0.05in}

$\bullet$ The sets $W(x,\omega )$ are {\em left-sequentially closed} in $(X,q)$ for all $(x,\omega )\in\gph\Omega $. To verify this property, it is sufficient to show that the limit of any convergent sequence $\{x_{k}\}\subset W(x,\omega)$ with $x_{k}\rightarrow x_{\ast}$ as $k\rightarrow\infty $ belongs to the set $W(x,\omega )$. By definition of $W$, find a sequence $\{\omega_{k}\}\subset\overline{\Omega}$ with $\omega_{k}\in\Omega(x_{k})$ for all $k\in\N$ satisfying
\begin{eqnarray*}
f(x_{k},\omega_{k})+q(x,x_{k})\xi\in f(x,\omega)-K\left[f\right].
\end{eqnarray*}
Since $\overline{\Omega}$ is a compact set, extract (without relabeling) a convergent subsequence from $\{\omega_{k}\}$ that converges to some
$\omega_{\ast}\in\overline{\Omega}$. This gives us $(x_{\ast},\omega_{\ast})\in\gph\Omega$ by the (left-sequential) closedness assumption on $\gph\Omega$. Then passing to limit with taking into account the level-closedness and lower semicontinuity assumptions imposed on $f$ and $q$ tells us that
\begin{eqnarray*}
f(x_{\ast},\omega_{\ast})+q(x,x_{\ast})\xi\in f(x,\omega)-K\left[f\right],\;\mbox{ i.e., }\;x_{\ast}\in W(x,\omega).
\end{eqnarray*}

$\bullet$ The sets $W(x,\omega )$ are {\em bounded from below} with respect to $\Theta+K\left[f\right]$ for all $(x,\omega)\in\gph\Omega$, where $f=f(x,\omega)$. Indeed, it follows from
\begin{eqnarray*}
W(x,\omega)\subset\left\{x^{\prime}\in X\;\mbox{ such that }\;q(x,x^{\prime})\xi\in f(x,\omega )-M-\Theta-K\left[f\right]\right\},
\end{eqnarray*}
where the bounded set $M$ is taken from the definition of the assumed quasiboundedness from below of the mapping $f$ with respect to the cone $\Theta $.\vspace*{0.05in}

$\bullet$ We have the inclusion $W(x^{\prime},\omega^{\prime})\subset W(x,\omega)$ for all $x^{\prime}\in W(x,\omega)$ and $\omega^{\prime}\in\overline{\Omega}$ with
\begin{eqnarray*}
f(x^{\prime},\omega^{\prime})+q(x,x^{\prime})\xi\le_{K\left[f\right]}f(x,\omega)
\end{eqnarray*}
To verify it, pick $x''\in W(x^{\prime},\omega^{\prime})$ and by construction of $W(x,\omega )$ in (\ref{0.4}) find $\omega''\in\Omega (x")$ satisfying the
inequality $f(x'',\omega'')+q(x^{\prime},x'')\xi\le_{K\left[f^{\prime}\right]}f(x^{\prime},\omega^{\prime})$. Summing up the last two inequalities and taking into account the triangle inequality for the quasimetric $q(x,x'')\le q(x,x^{\prime})+q(x^{\prime},x'')$, the choice of $\xi\in\Theta_{K}$ as well as the transitivity and convexity properties of $K$ ensuring that $K\left[f\right]+K\left[f^\prime\right]+\Theta_{K}=K\left[f\right]$, we get the relationships
\[
\begin{array}{lll}
&&f(x'',\omega'')+q(x,x'')\xi\\&=&\big(f(x^{\prime},\omega^{\prime})+q(x,x^{\prime})\xi\big)+\big(f(x'',\omega'')+q(x^{\prime},x'')\xi\big)\\
&&\hspace*{.01in}+\big(q(x,x'')-q(x,x^{\prime})-q(x^{\prime},x'')\big)\xi-f(x^{\prime},\omega^{\prime})\\
&\in&f(x,\omega)-K\left[f\right]+f(x^{\prime},\omega^{\prime})-K\left[f^{\prime}\right]-\Theta_{K}-f(x^{\prime},\omega^{\prime})\\
&\subset& f(x,\omega)-K\left[f\right].
\end{array}
\]
This clearly implies that $f(x'',\omega'')+q(x,x'')\xi\le_{K\left[f\right]}f(x,\omega)$, i.e., $x''\in W(x,\omega)$. Since $x''$ was chosen arbitrarily
in $W(x^{\prime},\omega^{\prime})$, we conclude that $W(x^{\prime},\omega^{\prime})\subset W(x,\omega)$.\vspace*{0.05in}

To proceed further, let us inductively construct a sequence of pairs $\{(x_{n},\omega_{n})\}\subset\gph\Omega$ and denote $f_{n}:=f(x_{n},\omega_{n})$ for
all $n\in\N\cup\{0\}$ by the following {\em iterative procedure}: starting with $(x_{0},\omega_{0})$ given in the theorem and having the $n$-iteration $(x_{n},\omega_{n})$, we select the next one $(x_{n+1},\omega_{n+1})$ by
\begin{eqnarray}\label{0.6}
\left\{\begin{array}{ll}
x_{n+1}\in W(x_{n},\omega_{n}),&\\
q(x_{n},x_{n+1})\disp\ge\sup_{x\in W(x_{n},\omega_{n})}q(x_{n},x)-\disp\frac{1}{n+1},&\\
\omega_{n+1}\in\Omega(x_{n+1}),\quad f(x_{n+1},\omega_{n+1})+q(x_{n},x_{n+1})\xi\;\le_{K\left[f_{n}\right]}f(x_{n},\omega_{n})\;&
\end{array}
\right.
\end{eqnarray}
It follows from construction (\ref{0.4}) of the sets $W(x,\omega)$ and their properties listed above that this iterative procedure is {\em well defined}. By
(\ref{0.5}) the sequence $\{f_{n}\}$ with $f_{n}:=f(x_{n},\omega_{n})$ is {\em nonincreasing} with respect to the ordering structure $K$ in the sense that $f_{n+1}\le_{K\left[f_{n}\right]}f_{n}$ for all $n\in\N\cup\left\{0\right\}$. Furthermore, the cone sequence $\{K(f_{n})\}$ is {\em nonexpansive}, i.e.,
\[
K\left[f_{n+1}\right]\subset K\left[f_{n}\right]\;\mbox{ for all }\;n\in\N\cup\left\{0\right\},
\]
which implies together with the convex-valuedness of $K$ that
\[
\sum_{n=0}^{m}K\left[f_{n}\right]=K\left[f_{0}\right]\;\mbox{ for all }\;m\in\N\cup\{0\}.
\]
Summing up the last inequality in (\ref{0.6}) from $n=0$ to $m$, we get with $t_{m}:=\sum_{n=0}^{m}q(x_{n},x_{n+1})$ that
\begin{eqnarray}\label{0.7}
t_{m}\xi\in f_{0}-f_{m+1}-K\left[f_{0}\right]\subset f_{0}-M-\Theta-K\left[f_{0}\right]\;\mbox{ for all }\;m\in\N\cup\{0\}.
\end{eqnarray}

Let us next prove by passing to the limit in \eqref{0.7} as $m\rightarrow\infty$ that
\begin{eqnarray}\label{0.8}
\sum_{n=0}^{\infty}q(x_{n},x_{n+1})<\infty.
\end{eqnarray}
Arguing by contradiction, suppose that (\ref{0.8}) does not hold, i.e., the increasing sequence $\left\{t_{m}\right\}$ tends to $\infty $ as $m\longrightarrow\infty$. By the first inclusion in (\ref{0.7}) and the boundedness of the set $M$ taken from the quasiboundedness of $f$ from below,
find a bounded sequence $\left\{w_{m}\right\}\subset f_{0}-M$ satisfying
\[
t_{m}\xi-w_{m}\in-\Theta-K\left[f_{0}\right],\;\mbox{ i.e., }\;\xi -w_{m}/t_{m}\in-\Theta-K\left[f_{0}\right],\quad m\in\N.
\]
Passing now to the limit as $m\longrightarrow\infty$ and taking into account the closedness of $\Theta$, the boundedness of $\left\{w_{m}\right\}$, and that $t_{m}\rightarrow\infty $, we arrive at $\xi\in-\Theta-K\left[ f_{0}\right]$ in contradiction to the choice of $\xi\in\Theta_{K}\setminus(-\Theta-K\left[ f_{0}\right])$. Thus (\ref{0.8}) holds and allows us for any $\varepsilon>0$ find a natural number $N_{\varepsilon}\in\N$ so that $t_{m}-t_{n}=\sum_{k=n}^{m-1}q(x_{k},x_{k+1})\le\varepsilon$ whenever $m\ge n\ge N_{\varepsilon}$. Hence
\[
q(x_{n},x_{m})\le\sum_{k=n}^{m-1}q(x_{k},x_{k+1})\le\varepsilon\;\mbox{ for all }\;m\ge n\ge N_{\varepsilon},
\]
which means that $\left\{x_{k}\right\}$ is a (left-sequential) {\em Cauchy sequence} in the quasimetric space $(X,q)$. Since $X$ is left-sequentially complete, there is $x_{\ast}\in X$ such that $x_{k}\longrightarrow x_{\ast}$ as $k\longrightarrow\infty$. Taking into account that $W(x_{k+1},\omega_{k+1})\subset W(x_{k},\omega_{k})$ and the choice of $x_{k+1}$, we get the estimate
\begin{eqnarray*}
\mbox{\rm radius\,}W(x_{k},\omega_{k}):=\sup_{x\in W(x_{k},\omega_{k})}q(x_{k},x)\le q(x_{k},x_{k+1})+\frac{1}{k+1}
\end{eqnarray*}
ensuring that $\mbox{\rm radius\,}W(x_{k},\omega _{k})\downarrow 0$ as $k\rightarrow\infty$. It follows from the left-sequential completeness of $X$ and the left-sequential closedness of $W(x_{k},\omega_{k})$ that
\begin{eqnarray}\label{0.9}
\bigcap_{k=0}^{\infty }W(x_{k},\omega_{k})=\left\{x_{\ast}\right\}\;\mbox{ for some}\;x_{\ast}\in X.
\end{eqnarray}

Now we justify the existence of $\omega_{\ast}\in\Omega(x_{\ast})$ such that $f_{\ast}:=f(x_{\ast},\omega _{\ast})\in\Min\big(F(x_{\ast}),K\left[f_{\ast}\right]\big)$ satisfies the relationships in (\ref{0.1}) and (\ref{0.2}). For each pair $(x_{k},\omega_{k})\in\gph\Omega$, define a subset of $\Bar{\Omega}$ by
\begin{eqnarray}\label{0.10}
R(x_{k},\omega_{k}):=\left\{\omega\in\Omega(x_{\ast})\big|\;f(x_{\ast},\omega)+q(x_{k},x_{\ast})\xi\le_{K\left[f_k \right]}f(x_{k},\omega_{k})\right\},\quad k\in\N.
\end{eqnarray}
Then we have the following properties:\vspace*{0.05in}

$\bullet$ The set $R(x_{k},\omega_{k})$ is {\em nonempty} and {\em closed} for any $k\in\N\cup\{0\}$ under the assumptions made. Indeed, the
nonemptiness follows directly from $x_{\ast}\in W(x_{k},\omega_{k})$ and the definition of $W$ in (\ref{0.4}). The closedness property holds since
$R(x_{k},\omega_{k})$ is the $(f(x_{k},\omega_{k})-q(x_{k},x_{\ast})\xi )$-level set of the mapping $f(x_{\ast},\cdot)$ with respect to the closed and convex cone $K\left[f_k\right]$ and since $f(x_\ast,\cdot)$ is assumed to be continuous. Furthermore, it follows from the inclusion $R(x_{k},\omega_{k})\subset\Omega(x_{\ast})\subset\overline{\Omega}$ and the compactness of $\Bar\O$ that $R(x_{k},\omega_{k})$ is a compact subset as well.\vspace*{0.05in}

$\bullet$ The sequence $\left\{R(x_{k},\omega_{k})\right\}$ is {\em nonexpansive}. To verify it, pick any $w\in R(x_{k+1},\omega_{k+1})$ and get
\[
f(x_{\ast},w)+q(x_{k+1},x_{\ast})\xi\in f(x_{k+1},\omega_{k+1})+K\left[f_{k+1}\right].
\]
Combining this with (\ref{0.6}) and then using the quasimetric triangle inequality together with the equality $K\left[f_{k+1}\right]+K\left[f_{k}\right]=K\left[f_{k}\right]$ tell us that
\begin{eqnarray*}
f(x_{\ast},w)+q(x_{k},x_{\ast})\xi\in f(x_{k},\omega_{k})+K\left[f_{k}\right],\;\mbox{ i.e., }\;w\in R(x_{k},\omega _{k}),
\end{eqnarray*}
which therefore justifies that $w\in R(x_{k+1},\omega_{k+1})\subset R(x_{k},\omega_{k})$.\vspace*{0.1in}

It follows from the properties of $R(\cdot,\cdot)$ established above and the compactness of $\overline{\Omega}$ that there exists $\overline{\omega}\in\overline{\Omega}$ satisfying the inclusion
\begin{eqnarray}\label{0.11}
\bar\omega\in\bigcap_{k=0}^{\infty}R(x_{k},\omega_{k}).
\end{eqnarray}
Denoting $f_{\bar\omega}:=f(x_{\ast},\bar\omega)$ and forming the $f_{\bar\omega}$--level set of $f(x_{\ast},\cdot)$ over $\Omega(x_{\ast})$ by
\begin{eqnarray}\label{0.12}
\Xi:=\left\{\omega\in\Omega(x_{\ast})\big|\;f_\omega:=f(x_{\ast},\omega)\le_{K\left[f_{\bar\omega}\right]}f(x_{\ast},
\bar\omega)=:f_{\bar\omega}\right\},
\end{eqnarray}
we obviously have that $\Xi$ is compact with $\bar\omega\in\Xi $. Employ now \cite[Corollary~5.10]{l89}, which ensures in our setting the existence of $\omega_{\ast}\in\Xi$ such that
\begin{eqnarray*}
f_{\ast }=f(x_{\ast},\omega_{\ast})\in\Min\big(f(x_{\ast },\Xi ),K\left[f_{\bar\omega}\right])\;\mbox{ with }\;f(x_{\ast},\Xi):=\bigcup_{\omega\in\Xi} \left\{f_{\omega}:=f(x_{\ast},\omega )\in P\right\}.
\end{eqnarray*}
This reads by the definition of Pareto efficiency that
\[
\big(f_{\ast}-K\left[f_{\bar\omega}\right]\big)\cap f(x_{\ast},\Xi)=\left\{f_{\ast}\right\}.
\]
Since $f_\ast\le_{K\left[f_{\bar\omega}\right]}f_{\bar\omega}$, we have $K\left[f_{\ast}\right]\subset K\left[f_{\bar\omega}\right]$; cf.\ the justifications for (\ref{0.5}). Thus we get
\[
(f_{\ast}-K\left[f_{\ast}\right])\cap f(x_{\ast},\Xi)=\left\{f_{\ast}\right\},\;\mbox{ i.e., }\;f_{\ast}\in\Min(f(x_{\ast },\Xi ),K\left[f_{\ast}\right]).
\]
Actually the following stronger conclusion holds:
\begin{eqnarray}\label{str}
f_{\ast}\in\Min(F(x_{\ast}),K\left[f_{\ast}\right])\;\mbox{ with }\;F(x_{\ast})=f\big(x_{\ast},\Omega(x_{\ast})\big)\supset f(x_{\ast},\Xi).
\end{eqnarray}
Arguing by contradiction, suppose that \eqref{str} does not hold and find $\omega\in\Omega(x_{\ast})\setminus\Xi$ such that $f_{\omega}\le_{K\left[f_{\ast}\right]}f_{\ast}$. Since $\omega_{\ast}\in\Xi$, we have $f_{\ast}\le_{K\left[f_{\bar\omega}\right]}f_{\bar\omega}$. Then the transitivity assumption {\bf(H3)} ensures that $f_{\omega}\le_{K\left[f_{\bar\omega}\right]}f_{\bar\omega}$, and so $\omega\in\Xi$ contradicting the choice of $\omega\in\Omega(x_{\ast})\setminus\Xi$. This justifies \eqref{str}.\vspace*{.05in}

Now we are ready to show that the pair $(x_{\ast},\omega_{\ast})$ satisfies the conclusions (\ref{0.1}) and (\ref{0.2}) of our variational principle. The inequality in (\ref{0.1}) immediately follows from $\omega_{\ast}\in R(x_{0},\omega_{0})$. To verify (\ref{0.2}), suppose the contrary and find a pair $(x,\omega)\in\gph\Omega$ with $f(x,\omega ) \neq f(x_{\ast},\omega_{\ast })$ satisfying
\begin{eqnarray}\label{0.13}
f(x,\omega)+q(x_{\ast},x)\xi \in f(x_{\ast},\omega_{\ast})+K\left[f_{\ast}\right].
\end{eqnarray}
Fix $k\in\N\cup\left\{0\right\}$ and sum up the three inequalities: (\ref{0.13}), (\ref{0.12}) with $\omega=\omega_{\ast}$, and (\ref{0.10}) with
$\omega=\bar\omega$. This gives us, by taking into account the triangle inequality as well as the relationships $f_{\ast}\le_{K\left[f_{\bar\omega}\right]}f_{\bar\omega}\le_{K\left[f_{k}\right]}f_{k}$, $K\left[f_{\ast}\right]\subset K\left[f_{\bar\omega}\right]\subset K\left[f_{k}\right]$, and $K\left[f_{\ast}\right]+K\left[f_{\omega}\right]+K \left[f_{k}\right]=K\left[f_{k}\right]$, that
\[
f(x,\omega)+q(x_{k},x)\xi\in f(x_{k},\omega_{k})-K\left[ f_{k}\right],\;\mbox{ i.e., }\;x\in W(x_{k},\omega_{k}),\quad k\in\N.
\]
This means that $x$ belongs to the set intersection in (\ref{0.9}), and thus $x=x_{\ast}$. Substituting it into (\ref{0.13}), we obviously get $f(x_{\ast},\omega)+q(x_{\ast},x_{\ast})\xi\in f(x_{\ast},\omega_{\ast})+K\left[f_{\ast}\right]$ and reduce it to
$$
(x_{\ast},\omega)\in f(x_{\ast},\omega_{\ast})-K\left[f_{\ast}\right],\;\mbox{ i.e., }\;f(x_{\ast},\omega)\le_{K\left[f_{\ast}\right]}f(x_{\ast},\omega_{\ast}).
$$
The latter shows that $f(x_{\ast},\omega)=f(x_{\ast},\omega_{\ast})\in\Min(F(x_{\ast}),K\left[f_{\ast}\right])$, which contradicts the assumption of $f(x,\omega)\ne f(x_{\ast},\omega_{\ast})$ and hence justifies \eqref{0.2}.\vspace*{0.05in}

To complete the proof of the theorem, it remains to estimate the distance $q(x_{0},x_{\ast})$ in \eqref{0.3} when $(x_{0},\omega_{0})$ is an approximate $\varepsilon\xi$--minimizer of $f$ over $\gph\Omega$. Arguing by contradiction, suppose that (\ref{0.3}) does not hold, i.e., $q(x_{0},x_{\ast})>\lambda$. Since $x_{\ast}\in W(x_{0},\omega_{0})$, we have
\begin{eqnarray*}
f(x_{\ast},\omega_{\ast})+\frac{\varepsilon}{\lambda}q(x_{\ast},x_{0})\xi\in f(x_{\ast},\omega_{\ast})-K\left[f_{\ast}\right],
\end{eqnarray*}
which together with $f(x_{\ast},\omega_{\ast})\le_{K\left[f_0\right]}f(x_0,\omega_0)$ yields $f(x_{\ast},\omega_{\ast})+\varepsilon\xi\in f(x_{\ast},\omega_{\ast})-K\left[f_{0}\right].$ This contradicts the approximate minimality of $(x_{0},\omega_{0})$ and thus ends the proof. $\h$\vspace*{0.05in}

Finally in this section, we present a direct consequence of Theorem~\ref{EVP-VOS} for the case when the mapping $f$ does not depend on the control variable $\omega$, which also provides a new variational principle for systems with variable ordering structures and is used in what follows.

\begin{Corollary}{\bf (variational principle for parameter-independent mappings with respect to variable ordering).}\label{EVP-Cor1} Let $f=f(x)$ be a mapping from $X$ to $P$ with $\dom f\ne\emp$ in the setting of Theorem~{\rm\ref{EVP-VOS}} under the assumptions made therein. Then for any $\varepsilon>0$, $\lambda>0$, $x_{0}\in\gph\Omega$, and $\xi\in\Theta_{K}\setminus(-\Theta-K\left[f(x_0)\right])$ with $\left\Vert\xi\right\Vert=1$ there is a point $x_{\ast}\in\dom f$ satisfying the relationships
\begin{eqnarray}\label{0.14}
f(x_{\ast})+\frac{\varepsilon}{\lambda}q(x_{0},x_{\ast})\xi\le_{K\left[f(x_{0})\right]}f(x_{0}),
\end{eqnarray}
\begin{eqnarray}\label{0.15}
f(x)+\frac{\varepsilon}{\lambda}q(x_{\ast},x)\xi\nleq_{K\left[f(x_{\ast})\right]}f(x_{\ast})\;\mbox{ for all }\;x\in\dom f\;\mbox{ with }\;f(x)\ne f(x_{\ast}).
\end{eqnarray}
If furthermore $x_{0}$ is an approximate $\varepsilon\xi$--minimizer of $f$ with respect to $K\left[f(x_0)\right]$, then $x_{\ast}$ can be chosen so that
in addition to \eqref{0.14} and \eqref{0.15} we have the estimate \eqref{0.3}.
\end{Corollary}
{\bf Proof.} It follows from Theorem~\ref{EVP-VOS} applied to the mapping $\Tilde{f}:X\times\overline{\Omega}\rightarrow P$ with $\Tilde{f}(x,\omega)=f(x)$ and a (compact) set $\overline{\Omega}$ consisting of just one point, say $\left\{\omega_{\ast}\right\}$. $\h$

\section{Applications to Goal Systems in Psychology}
\subsection{What Our Variational Principle Add to Goal System Theory}
{\bf (A) Variational rationality via variational analysis.} In the context of the variational rationality framework of \cite{s09,s10}, the new variational principle of Theorem~\ref{EVP-VOS} shows that, considering an {\em adaptive goal system} endowed with {\em variable cone-valued preferences} in the payoff space and a {\em quasimetric} on the space of means under (fairly natural) hypotheses of the theorem and starting from any feasible ``means-way of using these means" pair $\phi_{0}=(x_{0},\omega_{0})\in\Phi$, there exists a {\em succession of worthwhile changes} $\phi_{n+1}\in W(\phi_{n})$ with $n\in\N$, which ends at some {\em variational trap} $\phi_{\ast}=(x_{\ast},\omega_{\ast})\in\Phi$, where the agent {\em prefers to stay than to move}. The meaning of this is as follows; see the notation and psychological description in Section~2.

{\bf (i) Reachability and acceptability aspects along the transition}: we have $\phi_{\ast}\in W(\phi_{0})$, i.e., it {\em is worthwhile} to move directly from the starting means-ends pair $\phi_{0}$ to the ending one $\phi_{\ast}$.

{\bf (ii) Stability aspect at the end}: we have $W(\phi_{\ast})=\left\{\phi_{\ast}\right\}\Longleftrightarrow\phi\notin W(\phi_{\ast})$ for any
$\phi\in\Phi$, $\phi\ne\phi_{\ast}$, meaning that it is {\em not worthwhile} to move from the means-ends pair $\phi_{\ast }$ to a different one.

{\bf (iii) Feasibility aspect along the transition}: if $\phi_{0}=(x_{0},\omega_{0})\in\Phi$ is any $\varepsilon\xi$--{\em approximate minimizer} to $G(\cdot)$, then $x_{\ast}$ can be chosen such that in addition to {\bf(i)} and {\bf (ii)} we have $C(x_{0},x_{\ast})\le\lambda$.

{\bf (iv) The end is efficient as a Pareto optimal solution}. This is shown in the proof of Theorem~\ref{EVP-VOS} and discussed above.\\[1ex]
{\bf (B) When proof says more than statement.} Analyzing the statement and the proof of our variational principle in Theorem~\ref{EVP-VOS}, we can observe that--besides the variational trap interpretations, which are discussed above in {\bf(A)} and follow from the {\em statement} of the theorem--the {\em proof} itself offers much more from the psychological point of view. Indeed, the statement of Theorem~\ref{EVP-VOS} is an {\em existence result} while the proof provides a constructive {\em dynamical process}, which leads us to a solution. From the mathematical viewpoint, the situation is similar to the classical Ekeland principle with the proof given in \cite{e79}. From the psychological viewpoint, this is in accordance with the message popularized by Simon \cite{s55}: {\em decision and making process matters and can determine the end}. It is also a major point of the variational rationality approach \cite{s09,s10}: {\em to explain human desirable ends requires to exhibit human behavioral processes that can lead to them}. It senses that desirable ends must be reachable in an acceptable way by using feasible means. In other words, if the agent starts from any ``mean-ways of using these means" pair, pursues his/her goals by exploring enough each step and performing a succession of worthwhile changes or stays, then he/she will end in a strong behavioral trap, i.e., a Pareto solution more preferable to stay than to move even without any resistance to change. The given proof of Theorem~\ref{EVP-VOS} reveals at least {\em four} very important points discussed below in the rest of this subsection.\\[1ex]
{\bf(C) Worthwhile to change processes}. Parallel to \cite{s09} with using \cite{e79} in the case of nonadaptive models, the proof of Theorem~\ref{EVP-VOS} for the {\em adaptive} psychological models under consideration shows how the agent explicitly forms at each step a ``consideration set" to evaluate and balance his/her current motivation and resistance to change ``exploring enough" within the current worthwhile to change set trying to ``improve it enough" by inductively constructing a sequence of feasible pairs $\{(x_{n},\omega_{n})\}\subset\gph\Omega$. This nicely fits the famous concept of ``consideration sets" in marketing sciences defined first as ``evoked sets" by Howard and Sheth \cite{hs69}. The idea is that, at any given consumption occasion, consumers do not consider all the brands available while the current consideration/relevant set represents ``those brands that the consumer considers seriously when making a purchase and/or
consumption decision" as discussed, e.g., in \cite{b80,bl81}. The size of the consideration set is usually {\em small} relative to the total number of brands, which the consumer could evaluate. Then, using various heuristics, the consumer tries to simplify his/her decision environment.\\[1ex]
{\bf(D) Variational traps as desirable ends}. The worthwhile to change dynamical process given in the proof of Theorem~\ref{EVP-VOS} allows the agent to reach a
{\em variational trap} by a succession of worthwhile changes. In this model, it is a {\em Pareto ``means-way of using these means" pair}. Such a variational trap is related to two important concepts, {\em aspiration points} and {\em efficient points}, at the individual and collective levels discussed as follows.

{\bf (a) Aspiration points and Pareto points}. The Pareto point achieved in Theorem~\ref{EVP-VOS} is an aspiration point as defined in \cite{s09,s10} and then further studied and applied in \cite{fls12,ls12}. An {\em aspiration point} is such that, starting from any point of the worthwhile to change process, it is
worthwhile to move directly (in an acceptable way) to the given point of aspiration. It represents the ``{\em rather easy to reach}" aspect of a variational trap
while the other one (``difficult to leave") is more traditional as an equilibrium or stability condition.

{\bf (b) Optimal solutions}. The proof developed in Theorem~\ref{EVP-VOS} allows to study other types of {\em optimal solutions/minimal points}; compare, e.g., \cite{bm07,bm10,grtz03,q12} for various notions of this kind employing iterative procedures to derive for them variational principles of the Ekeland type.

{\bf (c) Individual or collective aspects: agents versus organizations}. In this paper we focus on the individual aspects of goal systems. The case of
organizations requires to consider {\em bilevel optimization problems} with leaders and followers. This will be a subject of our future research.\\[1ex]
{\bf (E) Variable preferences and efficiency for course pursuit processes}. Variable preferences can take different forms; we discuss it in more details in Appendix~4. In this paper we pay the main attention to ordering structures defined by {\em variable cone-valued preferences}. Among other types of variable preferences we mention {\em attention based preferences} discussed in \cite{b13,bg12}. Such variable preferences can be modelized and resolved by the approach developed in the proof of Theorem~\ref{EVP-VOS}.\\[1ex]
{\bf (F) Habituation processes as ends of course pursuit problems.} Our results help to modelize the emergence of {\em multiobjective habituation processes} with variable preferences. Such a formulation can represent agents who follow a habituation process with {\em multiple goals} as well as an organization, where each agent can have different goals. Then the procedure developed in the proof of Theorem~\ref{EVP-VOS} ends in a variational trap, which is a goal system habit for agents or a bundle of routines for organizations. This represents a habituation process in various areas of life, which is characterized by several properties such as repetitions, automaticity, control and economizing, etc.; see \cite{b94,s09,s10} for more more details and discussions.

\subsection{When Costs to Be Able to Change ``Ways of Using Means" Do Matter}

Consider a more general problem to change from a ``means-way of using these means" feasible pair $\phi=(x,\omega)$ with $\omega\in\Omega(x)\subset\overline{\Omega }$ to a new pair $\phi^{\prime}=(x^{\prime},\omega^{\prime})$ with $\omega^{\prime}\in\Omega(x^{\prime})\subset\overline{\Omega}$. In this general behavioral case, the full costs $C\left[(x,\omega),(x^{\prime},\omega^{\prime})\right]=C\left[\phi,\phi^{\prime}\right]$ to be able to change from a feasible pair $\phi=(x,\omega)$ with $\omega\in\Omega(x)$ to another feasible pair $\phi^{\prime}=(x^{\prime},\omega^{\prime})$ with $\omega^{\prime}\in\Omega(x^{\prime })$ must include the {\em two kinds of costs} in the sum: $C\left[(x,\omega),(x^{\prime},\omega^{\prime})\right]=C_{X}(x,x^{\prime})+C_{\Omega}(\omega,\omega^{\prime})\in P$.

Suppose now in the line of Theorem~\ref{EVP-VOS} that such vectorial costs are {\em proportional} to a vector $\xi\in P$, i.e., $C\left[(x,\omega),(x^{\prime},\omega ^{\prime})\right]=q\left[\phi,\phi^{\prime}\right]\xi$, where $q\left[\phi,\phi^{\prime}\right]\in\R_{+}$ is a quasidistance on $\overline{\Phi}:=X\times\overline{\Omega}$. This quasidistance modelizes the {\em total costs} to be able to change from one pair to another.

Let $\Phi:=\left\{(x,\omega),\omega\in\Omega (x)\right\}\subset\overline{\Phi }$ be the subset of feasible ``means-way of using these means" pairs. Then the {\em worthwhile to changes preference} over all the ``means-way of using these means" pairs $\phi=(x,\omega)\in\Phi$ is
\begin{eqnarray*}
\phi^{\prime}\ge_{K\left[f(\phi)\right]}\phi&\Longleftrightarrow&(x^{\prime},\omega^{\prime})\ge_{K\left[f(x,\omega)\right]}(x,\omega)\\
&\Longleftrightarrow&f(x^{\prime},\omega^{\prime})+C\left[(x,\omega),(x^{\prime},\omega^{\prime})\right]\le_{K\left[f(x,\omega)\right]}f(x,\omega)\\
&\Longleftrightarrow&f(\phi^{\prime})+C\left[\phi,\phi^{\prime}\right]\le_{K\left[f(\phi)\right]}f(\phi),
\end{eqnarray*}
where $C\left[(x,\omega),(x^{\prime},\omega^{\prime})\right]=q\left[(x,\omega),(x^{\prime},\omega^{\prime})\right]\xi$ while the pairs $\phi=(x,\omega)\in\Phi \subset\overline{\Phi}$ and $\phi^{\prime }=(x^{\prime},\omega^{\prime})\in\Phi\subset\overline{\Phi}$ are feasible, i.e., $\omega\in\Omega(x)\subset\overline{\Omega}$ and $\omega^{\prime}\in\Omega(x^{\prime})\subset\overline{\Omega}$. In this general case, the previous worthwhile to change sets read as follows:
\begin{eqnarray*}
W(x,\omega)&=&\left\{x^{\prime}\in X\big|\;\exists\;\omega^{\prime}\in\Omega(x^{\prime})\;\mbox{ with }\;
(x^{\prime},\omega^{\prime})\ge_{K\left[f(x,\omega)\right]}(x,\omega)\right\}\\
&=&\left\{x^{\prime}\in X\big|\;\exists\;\omega^{\prime}\in\Omega(x^{\prime})\;\mbox{ with }\;f(x^{\prime},\omega^{\prime})+C\left[(x,\omega),(x^{\prime},\omega^{\prime})\right]\le_{K\left[f(x,\omega)\right]}f(x,\omega )\right\}.
\end{eqnarray*}
Instead, we consider now the {\em new worthwhile to change set} defined by
\begin{eqnarray*}
W(\phi):&=&\left\{\phi^{\prime}\in\Phi\big|\;\phi^{\prime}\ge_{K\left[f(\phi)\right]}\phi\right\}=\left\{\phi^{\prime}\in\Phi\big|\;f(\phi^{\prime})+q\left[\phi,
\phi^{\prime}\right]\xi\le_{K\left[f(\phi)\right]}f(\phi)\right\}\\
&=&\left\{(x^{\prime},\omega^{\prime})\in\Phi\big|\;f(x^{\prime},\omega^{\prime})+q\left[(x,\omega),(x^{\prime},\omega^{\prime})\right]\xi
\le_{K\left[f(x,\omega)\right]}f(x,\omega)\right\}.
\end{eqnarray*}

In this setting, we can apply Corollary~\ref{EVP-Cor1}, where we replace the means $x\in X$ by the ``means-way of using these means" pairs $\phi= (x,\omega)\in\overline{\Phi}\subset X\times\overline{\Omega}$ to modelize such a situation. This variant has the following {\em two advantages}: {\bf(i)} it helps to modelize goal systems, where changing the ways of using means is costly; {\bf(ii)} it allows us to {\em drop the compactness} assumption on the set $\Bar\O$ of ways as in Theorem~\ref{EVP-VOS}. Now the state space is that of pairs $\phi=(x,\omega)\in\overline{\Phi}$, the vectorial payoff mapping is that of
unsatisfied needs $f:\phi\in\overline{\Phi}\longmapsto f(\phi)\in P$, and the real function $q(\phi,\phi^{\prime})\in\R_{+}$ denoted the quasidistance
between two pairs of ``means-way of using these means." Then, in this context of ``means-way of using means" pairs, we reformulate Corollary~\ref{EVP-Cor1}
as follows.
\begin{Corollary}{\bf (variational principle in ``means-way of using these means" setting).}\label{cor2} Let $(\overline{\Phi},q)$ be a left-sequentially complete quasimetric space, and let $K:P\rightrightarrows P$ be a cone-valued ordering structure satisfying assumptions {\bf(H2)} and {\bf(H3)}. Consider a mapping
$f:\Bar\Phi\rightarrow P$ with $\dom f\ne\emp$ being a left-sequentially closed subset of $\overline{\Phi}$. Assume also that:

{\bf(A1)} $f$ is quasibounded from below on $\dom f$ with respect to a convex cone $\Theta$.

{\bf(A2)} $f(\cdot)+\frac{\varepsilon}{\lambda}q(\phi,\cdot)$ is level-closed with respect to $K(\cdot)$ for all $\phi\in\overline{\Phi}$ and $\varepsilon,\lambda>0$.\\[1ex]
Then for any $\varepsilon$, $\lm>0$, $\phi_{0}\in\dom f$, and $\xi\in\Theta_{K}\setminus(-\Theta-K\left[f_{0})\right])$ with $\left\Vert\xi\right\Vert=1$ and $f_{0}:=f(\phi_{0})$ there is a point $\phi_{\ast}\in\dom f$,
satisfying the relationships
\[
f(\phi_{\ast})+\frac{\varepsilon}{\lambda}q(\phi_{0},\phi_{\ast})\xi\leq_{K\left[f_{0}\right]}f(\phi_{0}),
\]
\[
f(\phi)+\frac{\varepsilon}{\lambda}q(\phi_{\ast},\phi)\xi\nleq_{K\left[f_{\ast}\right]}f(\phi_{\ast})\;\mbox{ for all }\;\phi\in\dom f\;\mbox{ with }\;f(\phi)\ne f(\phi_{\ast}),
\]
where $f_{\ast}:=f(\phi_{\ast})$. If furthermore $\phi_{0}$ is an approximate $\varepsilon\xi$--minimizer of $f$ with respect to $K\left[f_{0}\right]$, then $\phi_{\ast}$ can be chosen so that in addition to \eqref{0.14} and \eqref{0.15} we have the estimate \eqref{0.3}.
\end{Corollary}
\textbf{Comments.} From the psychological point of view, Corollary~\ref{cor2} can be interpreted as follows. Starting from any ``means-way of using these means" pair $\phi_{0}\in\overline{\Phi}$, the agent who manages several goals by enduring costs to be able to change both the means used and the way of using them and whose next preference over the relative importance of each goal changes with the current pair $\phi$, can reach, in {\em only one worthwhile to change step}, a certain {\em variational trap} $\phi_{\ast}$, where it is not worthwhile to move. Moreover, given the desirability level $\varepsilon>0$ of the initial pair $\phi_{0}\in\overline{\Phi}$ and the size $\lambda>0$ of the limited resource, the agent accomplishes this worthwhile change in a {\em feasible way}, since the costs to be able to change $q(\phi_{\ast},\phi)$ are lower than the resource constraint $\lambda$.

\section{Conclusion}
The main mathematical result of this paper, Theorem~\ref{EVP-VOS}, as well as its consequences provide a far-going extension of the Ekeland variational principle aimed, first of all, to cover multiobjective problems with variable ordering structures. This major feature allows us to obtain new applications to adaptive psychological models within the variational rationality approach. Following this way, we plan to develop in our future research further applications of variational analysis to qualitative and algorithmic aspects of adaptive modeling in behavior sciences. One of our major attention in this respect is to extend the variational rationality approach and the corresponding tools variational analysis to decision making problems, where ``all things can be changed", i.e., with changeable decision sets, payoffs, goals, preferences, and contexts/parameters. Note that in a different setting, where decision sets and parameters can change along some  Markov chain, another approach to similar issues has been developed in the context of {\em habitual domain theory}; see \cite{ly09,ly11,ly12,yc10,yl09}.
A detailed comparison between the variational rationality approach and that of habitual domain theory has been recently given in \cite{bs13}.

\small\rm

\section{Appendix~1. Practical Means-Ends Rationality}
{\bf (A) Substantive rationality in economics and mathematics.} Simon \cite{s55} defines {\em rationality} as an {\em adequation between preestablished ends and some means to reach them}. A behavior (action or sequence of actions) is substantively rational when it allows the realization of some given desirable ends subject to given conditions and constraints. {\em Substantive rationality} represents perfect or global rationality, optimization evaluating the fit between desirable ends and feasible means in a comprehensive way. This includes the following {\em three steps}: {\bf(i)} the listing of all the possible alternatives/actions; {\bf(ii)} the determination of all the consequences that occur if the agent plans to adopt each of these alternatives in a deterministic or probabilistic way; {\bf(iii)} the evaluation of the consequences when the agent adopts these alternatives according to the preestablished ends, e.g., payoff functions, like costs, revenues, profits, and utilities. Then the agent is able to specify all the ends, to weight all of them, to examine and evaluate all the possible sets of means, to evaluate how well each given set of means achieves each end, to have the ability and resources to perform these evaluations, and finally to choose the set of means with the highest weighted score.\\[1ex]
{\bf(B) Incremental rationality: ``muddling through" in political sciences.} In administrative sciences,  Lindblom \cite{l59-1} considers incremental rationality via a non-global analysis. In this case we have:

$\bullet$ Ends and means are determined simultaneously since the agent knows his/her ends by considering the means, which the agent has in mind.

$\bullet$ To save the agent's limited resources (e.g., time, energy, money), many consequences are ignored since a full analysis is too costly. Then the evaluation of only major consequences should be provided. In the same vein, only a few means should be considered. Furthermore, evaluation of each of the considered means is incomplete: the consideration is only ``serious enough."

$\bullet$ The agent departs at each step from the status quo while not too much from it; hence the name ``branch" or ``muddling through" method. This means that he/she does not compare two new alternatives, but compares a new alternative to the old one (status quo) making small steps. Thus the agent, instead of using a rational comprehensive method, makes a finite succession of limited comparisons with respect to the status quo. In this way complex problems and decisions are significantly simplified.

$\bullet$ The choice among the means is determined by the agreement among interested parties when it concerns a group of agents. The focus is on {\em incremental objectives} while social objectives may have different value weights in different circumstances. Individuals may be unable to rank their own values when they are conflicting. Participants can also disagree on weights of critical values and even on sub-objectives.\\[1ex]
{\bf(C) Bounded and procedural rationality in decision theory.} In decision theory, Simon \cite{s55-1} considers {\em choice problems} instead of production problems, where agents, being bounded rational (they have limited resources and information), choose to be {\em procedural rational}. They search until these goals are satisfied by using heuristics for practical reasons and ignoring some alternatives to focus attention on a smaller subset of potential promising ones. In this way agents try to simplify their choice problems to economize their limited cognitive resources by taking into account their incomplete and inaccurate knowledge about the consequences of their actions. To simplify their choice, they accept to just {\em satisfy}, instead of to {\em optimize}, trying to find a course of action that is ``good enough," instead of being the best one. Since the goal to satisfy is less demanding, this relaxation procedure limits the sequential search for satisficing alternatives. Global/substantive rationality is only concerned with what is the result of the choice. Procedural rationality focuses on how the choice is done via a sequential search process, which stops when a satisfactory alternative is found. In a complex situation (e.g., too many alternatives and a long list of criteria required for a  sufficiently precise evaluation), the process the agent chooses--to be able at the second stage selecting a solution among so many alternatives--will change the final decision. Hence the final choice of a suitable alternative depends of the choice process, which uses algorithms, procedures and computations. A major goal in \cite{s55-1} is to discover the {\em symbolic processes} that people use in thinking, using an analogy between the computer and the human mind. Bounded rationality and procedural rationality are seen as complementarities. It is written in \cite{s55-1} that ``...bounded rationality does the critical part while procedural rationality does the assertive one..."\\[1ex]
{\bf (D) Problem solving in cognitive psychology.} In this paragraph, we take benefit of an excellent survey in Wikipedia (see ``Cognitive psychology and cognitive neuroscience. Problem solving from an evolutionary perspective"), which considers problem solving, problem finding, and problem shaping as  interrelated processes. Let us focus our attention to problem solving, which deals with any given situation that differs from a desired goal. This concerns the following major issues.

$\bullet$ {\bf Problem solving concept.} It is written in Wikipedia about this concept: ``Every problem is composed of an initial state, intermediate states, and a goal state (also: desired or final state) while the initial and goal states characterize the situations before and after solving the problem. The intermediate states describe any possible situation between initial and goal state. The set of operators builds up the transitions between the states. A solution is defined as the sequence of operators which leads from the initial state across intermediate states to the goal state."

$\bullet$ {\bf Difficult problems.} Again from there: ``Difficult problems have some typical characteristics that can be summarized as follows: intransparency (lack of clarity of the situation), commencement opacity, continuation opacity, polytely (multiple goals), inexpressiveness, opposition, transience, complexity (large numbers of items, interrelations and decisions), innumerability, connectivity (hierarchy relation, communication relation, allocation relation), heterogeneity, dynamics difficulties (time considerations), temporal constraints, temporal sensitivity, phase effects, dynamic unpredictability...The resolution of difficult problems requires a direct attack on each of these characteristics that are encountered."

$\bullet$ {\bf Well defined and ill defined problems.} Wikipedia says: ``{\em Well-defined} problems are such that it is possible to find an algorithmic solution. They can be properly formalized as: {\bf(i)} problems having clearly defined given states; {\bf (ii)} there is a finite set of operators, i.e., rules the agent may apply to given states; and {\bf (iii)} problems having clear goal states. For {\em ill-defined} problems (involving creativity) it is not possible to clearly define a given state and a goal state. Nevertheless, they often involve sub-problems that can be totally well-defined. Gestalt psychologists considered problem solving in situations requiring some novel means of attaining goals. In this context, problem solving requires a representation in a person's mind, and a reorganization or restructuring of this representation. Problem representations mean to model the situation as experienced by the agent to analyze it and split it into separate components: objects, predicates, state space, operators, selection criteria. In a goal-oriented situation, either the agent reproduces the response to the given problem from past experience (reproductive thinking), or the agent needs something new and different (insight) to achieve the goal. In this case, prior learning is of little help (productive thinking). Sometimes, previous experience or familiarity can even make problem solving more difficult. This is the case whenever habitual directions get in the way of finding new directions--an effect called fixation. ``...Functional fixedness concerns the solution of object-use problems. The basic idea is that when the usual way of using an object is emphasized, it will be far more difficult for a person to use that object in a novel manner."  Also ``...mental fixedness represents a person's tendency to respond to a given task in a manner based on past experience (mental set)..." This list goes on and on.

$\bullet$ {\bf Problem-solving strategies}. The aforementioned  Wikipedia paper details a long list of methods to solve a problem. ``The simplest method is to search for a solution by just trying one possibility after another." This method is often called {\em trial and error}. Means-end analysis provides yet another approach. It ``reduces the difference between initial state and goal state by creating subgoals until a subgoal can be reached directly." Problem-solving strategies include analogies, which ``describe similar structures and interconnect them to clarify and explain certain relations..."

Furthermore, ``...problem-solving strategies are the steps that one would use to find the problem(s) that are in the way to getting to one's own goal... In this
cycle one will recognize the problem, define the problem, develop a strategy to fix the problem, organize the knowledge of the problem, figure-out the resources at the user's disposal, monitor one's progress, and evaluate the solution for accuracy. Although called a cycle, one does not have to do each step in order to fix the problem, in fact those who don't are usually better at problem solving. The reason it is called a cycle is that once one is completed with a problem another usually will pop up...The following techniques are usually called problem-solving strategies: abstraction, analogy, brainstorming, divide and conquer: breaking down a large complex problem into smaller solvable problems, hypothesis testing, lateral thinking, means-ends analysis, method of focal objects (synthesizing seemingly non-matching characteristics of different objects into something new), morphological analysis (assessing the output and interactions of an entire system), proof (try to prove that the problem cannot be solved; the point where the proof fails will be the starting point for solving it), reduction (transforming the problem into another problem for which solutions exist), research (employing existing ideas or adapting existing solutions to similar problems), root cause analysis (identifying the cause of a problem), trial-and-error (testing possible solutions until the right one is found)..."\\[1ex]
{\bf(E) Practical rationality in philosophy and artificial intelligence.} In philosophy, practical rationality is the use of reasons to decide {\em how to
act}; namely, whether a prospective course of action is worth pursuing. {\em Practical reasoning} is the reasoning directed towards actions, i.e., deciding what to do; see Bratman \cite{b87-1}. It weighs conflicting considerations for and against competing options relative to the agent's desires, values, and believes.
{\em Theoretical/speculative reasoning} is the use of reasons to decide what to believe: the truth of contingent events as well as necessary truths. Finally, {\em productive/technical reasoning} attempts to find the best means for a given end. Means-end rationality advocates that agents have certain interests and rationality consists of acting to promote these interests by using some means. {\em Means-ends reasoning} is concerned with finding the means for achieving goals. The problem solver begins by focusing on the end/final goal and then determines a plan, or a strategy, for attaining the goal starting from his/her current situation. {\em Plan construction} is at the core of this backward consideration process. This method is used in artificial intelligence. It sets up backward smaller sub-goals, which complement the goal and then constantly re-evaluate the performance of those subgoals. By completing them, the agent approaches step by step the final goal. The overall goal is splitting into objectives, which in turn are splitting into individual steps or actions by taking into account that ``every attainable end is in itself a means to a more general end" (see Pollock \cite{p02-1}).

In artificial intelligence the belief-desire-intention models of agency (Bratman \cite{b87-1}, Rao and Georgeff \cite{rg95-1}, Georgeff et al. \cite{gptw99-1}) modelize practical reasoning as the succession of two main activities: {\bf (i)} {\em deliberation}, i.e., deciding what to do and {\bf (ii)} {\em means-ends reasoning}, i.e., deciding how to do it. Deliberation examines what an agent wants to achieve considering preferences, choosing goals, etc. Then deliberation generates intentions, i.e., plans that are turned into actions like the interface between deliberation and means-ends reasoning. Means-ends reasoning is used to
determine how the goals are to be achieved, e.g., thinking about suitable actions, resources, and how to organize activity.\\[1ex]
{\bf(F) Means-ends control beliefs in psychology.} In the context of the expectancy-valence theory of self regulation, Skinner et al. \cite{scb88-1} define the concept of {\em perceived control} in terms of three independent sets of beliefs, namely: {\em control beliefs} (expectations about the extent to which agents can obtain desired outcomes), {\em means-ends beliefs} (expectations about the extent to which certain potential causes produce
outcomes), and {\em agency beliefs} (expectations about the extent to which agents possess potential means).

\small\rm

\section{Appendix~2. Literature on Goal Systems}
{\bf (A) What is goal.} A goal is a mental representation of a desired future end, ``wished-for end that is considered to be attainable" (see Aarts and Elliot \cite{ae12-2}). An agent may have a lot of different goals, e.g., {\em approach or avoidance goals}, conscious or unconscious, distal (unrealistic, unattainable in the short run, abstract, vague, and long term like desires, aspirations, wishes, ideals) and {\em proximal goals} (realistic, attainable in the short run, precise, short term, measurable like a performance to reach, an intermediate means to get in order to satisfy a final goal), difficult or easy, etc. The content of a goal can be anything: an object, a person, a process, a way of doing, a means, a performance, or a learning goal, etc. A goal can be {\em individual} or {\em collective}. People set goals to motivate themselves toward desired ends. Goal can be related to {\em intrinsic motivations} or {\em extrinsic motivations}, achievements, ideals, oughts, etc. Goals are fundamental to help agents to self regulate their behaviors in order to be able to satisfy some of their needs and desires, to do what they want and need to do. {\em Self regulation} includes goal setting, goal striving, goal revision, and goal pursuit. Goal setting balances the degree of desirability and the degree of feasibility of a goal that must be both enjoyable (directly or indirectly relevant, rewarding) and realistic (specific, measurable, attainable and time bound). {\em Goal striving} describes the energy and persistent efforts an agent must spend to attain a goal along the path he follows to reach this goal (obstacles to be overcome, resistances to change, distractions, temptations, fatigue, discouragement, failed motivation, impatience, etc.). {\em Goal revision} concerns the ex post comparison (feedback) between the initial goal and the realized goal, the emotions this comparison can generate (deception, pain in case of failure, contentment, pleasure in case of success) and the choice of a new goal. People recalibrate the original goal and adjust in the wished-for end. ``Sometimes the shift is upward... Such an adjustment reflects a reaction to an insufficient challenge; a wished-for end that was previously not attainable has now been reached and is therefore supplanted by a higher level of aspiration... In contrast, sometimes the shift may be downward...Such an adjustment indicates a goal beyond reach and results in a less demanding target. Either way, many common activities involve similar dynamic revisions of goals across multiple time periods"
(see Wang and Mukhopadhyay \cite{wm12-2}). Then a new goal pursuit follows.\\[1ex]
{\bf (B) Goal structures.} In {\em management sciences}, Gutman was probably the first to consider means-ends chains as goal hierarchies; see \cite{g97-2} and the references therein. A {\em means-ends chain} (MEC) is defined as a {\em hierarchy} of goals that represents potential identities of the actions necessary for the person to reach his or her goal. Goals as ends in MECs can be grouped into three levels: {\em action goals} (concerned with the act itself), {\em outcome goals} (immediate effects of actions), and {\em consequences} (indirect effects stemming from outcomes). Then Pieters et al. \cite{pba95-2} examine a means-ends chain approach to goal-directed consumer behavior in terms of a hierarchical structure of increasingly more abstract goals, which are connected to one another through means-ends relationships. The goal structure incorporates both the relatively concrete level of specific action plans, which is concerned with the how of behavior, and the more abstract level of values and motives, which provides the ultimate reasons for pursuing a course of action and thus constitutes the way of
behavior. Vallacher and Wegner \cite{vw85-2} define a theory of action identification, where any action can be identified in many ways ranging from low-level identities that specify how the action is performed to high-level identities that signify why or with what effect the action is performed.\\[1ex]
{\bf(C) Goal systems.} In {\em psychology}, in the context of self regulation, Bandura \cite{b88-2} examines {\em hierarchical structures} of goal systems in relation to attribution theory, expectancy value theory, and goal theory, where self regulatory mechanisms, negative feedbacks, goal properties, self motivation, and affective consequences of goal discrepancies play a major role. Later in \cite{ksffcs02-2}, Kruglanski et al. define goal systems as mental representations of {\em motivational networks} composed of interconnected goals and means. As {\em mental representations}, these systems possess {\em cognitive properties}. As {\em desirable end and feasible means} structures, they have also {\em motivational properties} because the strength of an agent motivation increases with the
degree of desirability and feasibility of each goal. A goal system includes several desirable and different goals, which collaborate with each other to
help attaining a global state of satisfaction and which also compete with each other to mobilize mental and other means and resources. Then the degrees of desirability and feasibility of each different goal interact via a certain internal goal system {\em coopetition} as well as {\em competition} between goals.\\[1ex]
{\bf (D) Cognitive properties of goal systems.}

$\bullet$ {\bf Interconnectedness: its forms and its strength.} Goal systems consist of mentally represented networks of goals and related means helping to attain them, where these means may be associated with other means. The links between goals and means (as mental representations of motivational constructs) can be facilitative as well as inhibitory. {\em Facilitative links} represent vertical connections between goals and their corresponding means. {\em Inhibitory links} connect lateral elements, competing goals, and competing means. Top level goals are cognitively connected to their various subgoals ``en route to that goal" that are finally connected to their own means of attainment. {\em Lateral links} connect different subgoals and different attainment means to those particular subgoals. The form of the network depends of goals linked to a given mean (the shape of the multifinality set) and of means linked to a given goal (the shape of the equifinality set). The strength of the cognitive associations between goals and means depend of the size of these sets. The larger they are, the weaker are these cognitive interconnections. Another determinant of the strength of associations are repeated pairing of goals and means. Subconscious impacts can be explained in such a way that focal goals are explicit and background goals are not consciously registered. Then goals can be subliminally primed. Goal systems are highly flexible and context-dependent. Their shape and the strength of some cognitive associations may vary in accordance with situational framing effects. Some
goals and means can be activated in certain contexts, but not in some other ones.

$\bullet$ {\bf The allocations properties of goal systems} play a fundamental role because goal pursuit is resource dependent and mental resources are limited.
The more resources are accorded to reach a given goal, the less mental resources are left to reach the remaining goals. Then different goals compete for scarce resources and the same resources can be allocated to different goals.\\[1ex]
{\bf (E) Motivational properties of goal systems.}

$\bullet$ {\bf Goal striving}. To attain specific desirable goals within a goal system spends limited resources. Then the agent cannot reach all his/her goals
at the same time and must alternate allocations of resources between promotion and prevention goals. This implies that the successful attainment of
one's goal that engenders positive feelings of pleasure or satisfaction, can be compensated, or not, by the failure to attain one's goals that engenders
negative feelings of displeasure and disappointment. Higgins \cite{h97-2} shows that the attainment of {\em promotion goals} generates feelings of happiness and
pride while their non-attainment generates feelings of sadness and dejection. On the other hand, the attainment of {\em prevention goals} engenders feelings of calm and relaxation while their non-attainment gives rise to tension and agitation.

$\bullet$ {\bf Goal commitment} defines the degree to which an individual is determined to pursue a goal. It increases with the value assigned to the
goal and its expectancy of attainment. The presence of multiple goals and limited common means pose the problem of when to abandon a current goal and
pursue other ones, or persist in the pursuit of these current goals.

$\bullet$ {\bf Means choice and substitution.} The choice of focal goals to satisfy soon determines the allocations of means. At the same time, available means determine the focal goals, which can be satisfied in the short run. The choice of means depends on how they impact the expectancy of goal attainment. Background
goals can constraint the current allocation of resources (to save energy, time constraints, distractions, etc). Some means can be substituted into
others when they fail to help in reaching the goal. New means can be generated if necessary.

\section{Appendix~3. Stability/Stays and Change Dynamics}

Let us cite, among many other approaches, the following ones in behavioral sciences related to {\em stability and change dynamics}.

$\bullet$ {\bf In biology:} ``punctuated equilibria" and the {\em tempo and mode of evolution} in the lines of the seminal contribution by Simpson \cite{si44-3} to the evolutionary synthesis; see also Gould and Eldredge \cite{ge77-3}.

$\bullet$ {\bf In economics:} these dual dynamics of ``stays and changes" have been the main grand vision of Schumpeter \cite{s02-3,s05-3}. He  says that {\em
Statics} (equilibrium economics) and {\em Dynamics} (evolutionary economics) are completely different fields, which concern not only different problems but
also different methods and different materials. They are not two chapters of one and the same theoretical building but two completely independent buildings. Only Statics has hitherto been somewhat satisfactorily worked up and we essentially only deal with it in Schumpeter's work. Dynamics that is still at its beginnings is a ``land of the future" as discussed in Andersen \cite{a06-3,a08-3}. Among other publications on variable preferences, belief formation,
habit formation and related topics, we mention the book by Becker \cite{b93-3} with the references therein.

$\bullet$ {\bf In management sciences:} adaptive models are often described and cited but {\em not formalized}. Let us mention, among many others, organizational dynamics of stability and change, routine change, routine formation and evolution in Nelson and Winter \cite{nw82-3}, ambidextry behaviors of exploration-exploitation in March \cite{m91-3}, organizational learning in Argyris \cite{a99-3}, models for diagnosing organizational behavior in Nadler and Tushman \cite{nt80-3}, models of buyer behavior in Howard and Sheth \cite{hs69-3}, dynamics of consideration set formation in Brisoux \cite{b80-3} and Brisoux and Laroche \cite{bl81-3}, stage models of consumer behavior in Engel et al. \cite{ebm95-3}, change management in Bullock and Batten \cite{bb85-3} Carnall \cite{c90-3} and Hiatt \cite{h06-3}, transition management  in Bridges \cite{b91-3} and Bridges and Mitchell \cite{bm02-3}, learning organizations in Senge \cite{s90-3}, leader-follower models of change in Kotter \cite{k95-3}, pace models of stability and change that can be radical/revolutionary, incremental/gradual, continuous or episodic, slow-rapid, etc.

$\bullet$ {\bf In psychology:} central approaches are: {\em proactive} and {\em reactive behaviors}, {\em tension reduction-tension production dynamics} (see, e.g., Maslow \cite{m54-3},  Murray \cite{m38-3}, and Bandura \cite{b99-3}), {\em self regulation problems} (goal setting, goal striving, goal
revision, and goal pursuit as in Bandura \cite{b99-3}, Carver and Scheier \cite{cs98-3}, and De Ridder and de Wit \cite{rw06-3}) as well as {\em human development and enaction models of stability and change}. This includes: the three-stage process of ``unfreezing change and freezing" in Lewin \cite{l47-3,l51-3}, the model of planned change in Lippitt et al. \cite{lww58-3}, the action phase model in Heckhausen \cite{h07-3} and Gollwitzer \cite{g90-3}, the transtheoretical model of change in Prochaska and DiClemente \cite{pd92-3,pd94-3}, stability and change models of attitudes  in McGuire et al. \cite{mla85-3} and Eagly and Chaiken \cite{ec95-3}, vocational choice processes in Gottfredson \cite{g81-3}, experiential learning in Kolb \cite{k84-3}, gap models of change and managing change in Beckhard and Harris \cite{bh87-3}, goal representation and goal systems in Kruglanski et al. \cite{ksffcs02-3}, hope theory in Snyder \cite{s94-3}, habit formation and habit-goal interface in Wood and Neal \cite{wn07-3}, etc.

$\bullet$ {\bf In medicine:} the {\em health belief model} developed by Janz and Becker \cite{jb84-3}.

$\bullet$ {\bf In philosophy and artificial intelligence:} BDI (Belief--Desire--Intention) models of intention formation in \cite{b87-3} and Wooldridge \cite{w02-3},  ``how we think" models in  Schoenfeld \cite{s11-3}, etc.

$\bullet$ {\bf In political sciences:} see muddling through processes in Lindblmom \cite{l59-3} and various other models of the dynamic of conflicts.

$\bullet$ {\bf In sociology:} cultural changes and group dynamics, see, e.g., Lewin \cite{l47-3,l51-3} and  Schein \cite{s99-3} among other publications in this direction.

$\bullet$ {\bf In game theory:} various models with interrelated agents, convergence to Nash equilibria, etc.; see Chen and Gazzale \cite{cg04-3}, Soares et al. \cite{sbcg13}, and the references therein.

\section{Appendix~4. How Preferences Change}

In this paper, we suppose that preferences over vectorial payoffs, $f(x,\omega)\in P$ or $g(x,\omega)\in P$, can change. They were modelized  above by {\em variable cone-valued structures} $K\left[f(x,\omega)\right]$ or $K\left[g(x,\omega)\right]$, i.e., via $\ge_{K\left[f(x,\omega)\right]}$ or $\ge_{K\left[g(x,\omega)\right]}$. These preferences change because these cones depend on the {\em status quo payoff} given either  by the current vector of
unsatisfied needs $f_{k}=f(x_{k},\omega_{k})\in P$, or by the current vector of gains $g_{k}=g(x_{k},\omega_{k})\in P.$ Then the status quo payoff changes as long as the agent prefers to change than to stay at the current status quo using a worthwhile change. If the agent has not reached a {\em variational trap} at which he/she prefers to stay than to move in a worthwhile way, then the agent prefers to change. This means that the agent enters in a course pursuit between the old preference $\ge_{K\left[f(x_{k},\omega_{k})\right]}$ or $\ge_{K\left[g(x_{k},\omega_{k})\right]}$ and chooses a new position $(x_{k+1},\omega_{k+1})$, where each new position defines a new preference. In turn, each new preference implies the choice of a new worthwhile position. Hence this process of the preferences change is {\em bounded rational} since worthwhile changes are improving, while it is {\em not an optimization process} in general. At each step, the agent who tries to improve, in a worthwhile way with respect to the current preference and the status quo, largely simplifies the process of comparison. He/she only compares a new alternative with the status quo while does not compare two new alternatives to understand which of them is better. This is one of the main characteristic of a ``muddling through" process in Lindblom \cite{l59-4}.

It is shown in the economic and psychological literature that real-life preferences {\em depend on reference points}; see, e.g., Kahneman and Tversky \cite{kt79-4} for risky choices and Tversky and Kahneman \cite{tk91-4} for riskiness choices. These reference points determine the position of choice outcomes on a value function determining their desirability. Then preferences change because reference points change along the process of decisions to act and actions. There is considerable empirical evidences showing that preference depend simultaneously of multiple reference points; we mention Golman and Loewestein \cite{gl11-4} among other references. There can be any stimulus, which ``other stimuli are seen in relation to" (see Rosch \cite{r75-4}; in particular, current states, feasibility conditions linked to current endowments (physical, physiological, psychological, cognitive, emotional, and motivational resources), or desirability conditions linked to multiple reference points (a survival state, a foregone alternative, a future aspiration state or aspiration level, an expectation, a social comparison, etc.). More generally, variable preferences depend on internal states of the agent and his/her environment (external state including actions of other agents). Preferences can also change (although less) with a stock of habits; see Abel \cite{a90-4}, Becker and Murphy \cite{bm88-4}, Carroll \cite{c00-4}, Ravn et al. \cite{rsu06-4}, Wendner \cite{w02-4}, Wood and Neal \cite{wn07-4}, etc.

For more discussions on how multiple reference points change as well as on related topics, see (among others) the following bibliographies:

$\bullet$ {\em Distal goals} in Bandura and Schunk \cite{bs81-4}, {\em dreams and aspirations in psychology} in Cooke et al. \cite{cjcnw04-4}, {\em consumer choice in management} and {\em desire in philosophy and artificial intelligence} via the BDI approach in  Bratman \cite{b87-4} and Woodridge \cite{w00-4,wo02-4}. These references represent the best that the reader can think about in a lot of different domains of interest, without even expecting to reach them. We also refer the reader to Appadurai \cite{a04-4} for {\em capacity to aspire} in development economics and to Helson \cite{h64-4} and Easterlin \cite{e05-4} for {\em adaptation level theory} and for {\em well being theory} in positive psychology.

$\bullet$ {\em Expectations} defined as expected payoffs or expected incomes, which can act as reference points; see Abeler et al. \cite{afgh11-4} with the bibliographies therein.

$\bullet$ {\em Hopes in psychology} (see Snyder \cite{s94-4}), where hope means that ``you can get here from there."

$\bullet$ {\em Proximal goals}  (Bandura and Schunk \cite{bs81-4} and Heath et al. \cite{hlw99-4}) and {\em aspirations} defined in economics and management as satisficing levels or threshold points, which separate failure from success; see Simon \cite{s55-4}. We warn the reader that the proximal notion of aspirations in economics and management differs from the similar notion in psychology. This causes a lot of confusions in the literature about aspirations; see, e.g.,
Golman and Loewenstein \cite{gl11-4} and aspiration adaptation theory in Selten \cite{s88-4}.

$\bullet$ {\em Decision framing and mental accounting}, which decompose references points in Thaler \cite{t99-4}, {\em selective attention} in Bhatia \cite{b13-4} and Bhatia and Golman \cite{bg12-4}, {\em prototype theory} in Rosch \cite{r75-4}, etc.

$\bullet$ {\em Cognitive reference points} in management; see consumer decision making theory in Klein and Oglethorpe \cite{ko87-4} and the complexity of their formation.

$\bullet$ {\em Constructing preferences} from memory in Weber and Johnson \cite{wj06-4}.

$\bullet$ {\em Social comparisons}, like payoffs relative to a group and social status; see Abel \cite{a90-4}.

$\bullet$ Finally, let us give a short list of recent references on vectorial and set-valued versions of the Ekeland variational principle,  quasimetric spaces, and variable domination structures. From the mathematical viewpoint, our paper considers {\em variable preferences} as a fundamental aspect of any {\em adaptive behavioral model}. We show that the variational principles in Bao and Mordukhovich \cite{bm07-4,bm10-4} and the methods of their proofs can be extended to include
variable cone-valued preferences. We also mention some other recent extensions of the EVP given by Guti\'{e}rrez et al. \cite{gjn08-4}, Liu and Ng \cite{lm11-4}, Khanh and Quy \cite{kq11-4}, Qiu \cite{q12-4}, and Bae and Cho \cite{bc13-4}. The papers by  Abdeljawad et al. \cite{akb09-4} and by Cho and Bae \cite{cb12-4} contain more discussions and references on quasimetric (cone metric) spaces. There are just few recent references on variable domination structures; see Bao and Mordukhovich \cite{bm13-4} Chen and Yang \cite{cy02-4}, Eichfelder \cite{e11-4,e13-4}, and Eichfelder and Ha \cite{eh13-4}. In Luc and Soubeyran \cite{ls13-4}, the reader can find discussions on reference-point dependent preferences. To conclude this short survey, we refer the reader to Bell Cruz and Bouza Allende \cite{bb13-4} for the new steepest descent-like variable order vector optimization problem and a variational method of its solution. Note that none of the papers considered above examines, as we do, variational principles in quasimetric spaces with variable cone-valued structures. None of them gives applications in behavioral sciences as it is done in our paper.


\begin{thebibliography}{99}

\bibitem{ah96} Alber, S., Heward, W.:  Twenty-five behavior trapsguaranteed to extend your students' academic and social skills,
{\em Interven. School Clinic} {\bf 31}, 285--289 (1996)

\bibitem{a04} Appadurai, A.: The capacity to aspire: culture and the terms of recognition, in {\em Culture and Public Action} (eds. V. Rao and
M. Walton), pp.\ 59--84, The World Bank Publishers (2004)

\bibitem{abm05} Attouch, H., Buttazzo, G., Michaille, G.: {\em Variational Analysis in Sobolev and BV Spaces}, SIAM Publications (2005)

\bibitem{as11} Attouch, H., Soubeyran, A.: Local search proximal algorithms as decision dynamics with costs to move, {\em Set-Valued Var. Anal.}
{\bf 19}, 157--177 (2011)

\bibitem{bw70} Baer, D., Wolf, M.: The entry into natural communities of reinforcement, in: {\em Control of Human Behavior} (eds. R. Ulrich,
T. Stachnick, J. Mabry), pp.\ 319--324, Glenview, IL, Scott Foresman (1970)

\bibitem{b94} Barg, J.: The four horsemen of automaticity: awareness, intention, efficiency, and control in social cognition, {\em Handbook of
Social Cognition} {\bf 1}, 1-40 (1994

\bibitem{bm07} Bao, T.Q., Mordukhovich, B.S.: Variational principles for set-valued mappings with applications to multiobjective optimization,
{\em Control Cyber.} {\bf 36}, 531--562 (2007)

\bibitem{bm10} Bao, T.Q., Mordukhovich, B.S.: Relative Pareto minimizers for multiobjective problems: existence and optimality conditions, {\em Math.
Program.} {\bf 122}, 301--347 (2010)

\bibitem{bm10a} Bao, T.Q., Mordukhovich, B.S.: Set-valued optimization in welfare economics, {\em Adv. Math. Econ.} {\bf 13}, 113-153 (2010)

\bibitem{bm13} Bao, T.Q., Mordukhovich, B.S.: Necessary nondomination conditions for set and vector optimization with variable structures, {\em J.  Optim.
Theory Appl.}, DOI 10.1007/s10957-013-0332-6 (2013)

\bibitem{b02} Baumeister, R.: Ego-depletion and self-control failure: an energy model of the self's executive function, {\em  Self Identity} {\bf 1},
129--136 (2002)

\bibitem{bh96} Baumeister, R., Heatherton, T.: Self regulation failure: an overview, {\em Psychol. Inquiry} {\bf 7}, 1--15 (1996)

\bibitem{bb13} Bello Cruz, J., Bouza Allende, G.: A steepest descent-like method for variable order vector optimization problems, {\em J. Optim.
Theory Appl.}, to appear (2013)

\bibitem{bs13} Bento, G, Soubeyran, A.: Generalized inexact proximal algorithms: Habit's formation with resistance to change, following
worthwhile changes, {\em J. Optim. Theory Appl.}, to appear (2013)

\bibitem{b13} Bhatia, S.: Associations and the sccumulation of preference, {\em Psychol. Rev.}, to appear( 2013)

\bibitem{bg12} Bhatia, S., Golman, R.: Attention and reference dependence, Department of Social \& Decision Sciences, Carnegie Mellon
University (2012)

\bibitem{bz05} Borwein, J.M., Zhu, Q.J.: {\em Techniques of Variational Analysis}, Springer (2005)

\bibitem{b09} Bridges, W.: {\em Managing Transitions: Making the Most of Change}, Nicholas Brealey Publishing, 3rd edition (2009)

\bibitem{b80} Brisoux, J.: {\em Le Ph\'{e}nom\`{e}ne des Ensembles \'{e}Voqu\'{e}s: Une \'{E}tude Empirique des Dimensions Contenu et
Taille}, Th\`{e}se de doctorat, Universit\'{e} Laval (1980)

\bibitem{bl81} Brisoux, J.E., Laroche, M.: Evoked set formation and composition: an empirical investigation under a routinized response
behavior situation, in: {\em Advances in Consumer Research} (ed. K.B. Monroe), pp.\ 357--361, Ann Arbor, MI, Association for Consumer Research (1981)

\bibitem{cy02} Chen, G.Y., Yang, X.Q.: Characterizations of variable domination structures via nonlinear scalarization, {\em  J. Optim. Theory
Appl.} {\bf 112}, 97--110 (2002)

\bibitem{coss13} Cruz Neto, J.X., Oliveira, P.R., Soares Jr., P.A., Soubeyran, A.: Learning how to play Nash, potential games and
alternating minimization method for structured nonconvex problems on Riemannian manifolds, {\em J. Convex Anal.} {\bf 20}, 395--438 (2013)

\bibitem{e11} Eichfelder, G.: Optimal elements in vector optimization with a variable ordering structure, {\em J. Optim. Theory Appl.} {\bf 151},
217--240 (2011)

\bibitem{eh13} G. Eichfelder and T.X.D. Ha,  Optimality conditions for vector optimization problems with variable ordering structures,
{\em Optimization}, to appear (2013)

\bibitem{e72} Ekeland, I.: Sur les probl\'emes variationnels, {\em C. R. Acad. Sci. Paris} {\bf 275}, 1057--1059 (1972)

\bibitem{e74} Ekeland, I.: On the variational principle, {\em J. Math. Anal. Appl.} {\bf 47}, 324--353 (1974)

\bibitem{e79} Ekeland, I.: Nonconvex minimization problems, {\em Bull. Amer. Math. Soc.} {\bf 1}, 432--467 (1979)

\bibitem{fls12} Flores-Bazan, F., Luc, D.T., Soubeyran, A.: Maximal elements under reference-dependent preferences with applications to
behavioral traps and games, {\em J. Optim. Theory Appl.} {\bf 155}, 883--901 (2012)

\bibitem{ft13} Farokhinia, A., Taslim, L.: On assymmetrc distance, {|em J. Anal. Num. Theo.} {\bf 1}, 11--14 (2013)

\bibitem{grtz03} G\"{o}pfert, A., Riahi, H., Tammer, C., Z\u{a}linescu, C.: {\em Variational Methods in Partially Ordered Spaces},
Springer (2003)

\bibitem{gjn12} Guti\'{e}rrez, C., Jim\'{e}nez, B., Novo, V.: A set-valued Ekeland's variational principle in vector optimization,
{\em SIAM J. Control Optim.} {\bf 47}, 883--903 (2008)

\bibitem{hkr98} Hammond, J.S., Keeney, R.L., Raiffa, H.: The hidden traps in decision making, {\em Harvard Business Review} (1998)

\bibitem{hm06} Heifetz, A., Minelli, E.: Aspiration traps, Discussion Paper 0610, Universit\'{a} degli Studi di Brescia (2006)

\bibitem{hs69} Howard, J.A., Sheth, J.N.: {\em The Theory of Buyer Behavior (Marketing)}, John Wiley \& Sons (1969)

\bibitem{kq11} Khanh, P.Q., Quy D.N.: On generalized Ekeland's variational principle and equivalent formulations for set-valued mappings,
{\em J. Global Optim.} {\bf 49}, 381--396 (2011)

\bibitem{ly09} Larbani, M., Yu, P.L.: Two-person second-order games, Part II: Restructuring operations to reach a win-win profile, {\em J. Optim. Theory
 Appl.} {\bf 141}, 641--659 (2009)

\bibitem{ly11} Larbani, M., Yu, P.L.: $n$-Person second-order games: A paradigm shift in game theory, {\em J. Optim. Theory Appl.} {\bf 149}, 447--473 (2011)

\bibitem{ly12} Larbani M., Yu P.L.: Decision making and optimization in changeable spaces, a new paradigm, {\em J. Optim. Theory Appl.} {\bf 155},
727--761 (2012)

\bibitem{lm93} Levinthal, D., March, J.: The myopia of learning, {\em Strateg. Manag. J.} {\bf 14}, 95--112 (1993)

\bibitem{ln11} Liu, C.G., Ng, K.F.: Ekeland's variational principle for set-valued functions, {\em SIAM J. Optim.} {\bf 21}, 41--56 (2011)

\bibitem{l89}  Luc, D.T.: {\em Theory of Vector Optimization}, Springer (1989)

\bibitem{ls12} Luc, D.T., Soubeyran, A.: Variable preference relations: existence of maximal elements, {\em J. Math. Econ.}, to appear (2013)

\bibitem{m91} March, J.G.: Exploration and exploitation in organizational learning, {\em Organiz. Science} {\bf 2}, 71--87 (1991)

\bibitem{m06} Mordukhovich, B.S.: {\em Variational Analysis and Generalized Differentiation, I: Basic Theory, II: Applications}, Springer (2006)

\bibitem{mos11} Moreno, F.G., Oliveira, P.R., Soubeyran, A.: A proximal algorithm with quasidistance. Application to Habit's Formation,
{\em Optimization} {\bf 61}, 1383--1403 (2011)

\bibitem{p93} Plous, C.: {\em The Psychology of Judgment and Decision Making}, Mcgraw-Hill Book Company (1993)

\bibitem{q12} Qiu, J.H.: On Ha's version of set-valued Ekeland's variational principle, {\em Acta Math. Sinica} (English Series) {\bf 28},
717--726 (2012)

\bibitem{r06} Ray, D.: Aspirations, poverty and economic change, in: {\em What We Have Learnt about Poverty}  (eds. A. Banerjee, R. Benabou,
D. Mookherjee), Oxford University Press (2006)

\bibitem{rw98} Rockafellar, R.T., Wets, J.-B.: {\em Variational Analysis}, Springer (1998)

\bibitem{s55} Simon, H.: A behavioral model of rational choice, {\em Quarterly J. Econom.} {\bf 69}, 99--188 (1955)

\bibitem{s09} Soubeyran, A. : Variational rationality, a theory of individual stability and change, worthwhile and ambidextry behaviors, Preprint,
GREQAM, Aix-Marseille University (2009)

\bibitem{s10} Soubeyran, A.: Variational rationality and the unsatisfied man: routines and the course pursuit between aspirations, capabilities and
beliefs, Preprint, GREQAM, Aix-Marseille University (2010)

\bibitem{s04} Stephen, F.: {\em The Power of Reinforcement}, SUNY Press (2004)

\bibitem{w74} Walras, L.: {\em El\'{e}ments d'\'{E}conomie Politique Pure}, Lausanne, Corbaz Publishers (1874)

\bibitem{yc10} Yu, P.L., Chen,Y.C.: Dynamic MCDM, habitual domains and competence set analysis for effective decision making in changeable spaces, in {\em Trends
in Multiple Criteria Decision Analysis}, (eds. M. Ehrgott, J.R.F. Salvatore Greco), Chapter~1, Springer (2010)

\bibitem{yl09} Yu, P.L., Larbani M.: Two-person second-order games, Part 1: Formulation and transition anatomy, {\em J. Optim. Theory Appl.} {\bf 141},
619--639 (2009)
\end{thebibliography}

\begin{thebibliography}{99}

\bibitem{b87-1} Bratman, M.: {\em Intentions, Plans, and Practical Reason}, Harvard University Press (1987)

\bibitem{gptw99-1} Georgeff, M., Pell, B., Pollack, M., Tambe, M., Wooldridge, M.: The belief-desire-intention model of agency, {\em Proc. Agent. Theor.
Architect. Lang.} (1999)

\bibitem{l59-1} Lindblom, C.: The science of ``muddling through", {\em Public Admin. Rev.}, 79--88 (1959)

\bibitem{p02-1} Pollock, J.: The logical foundations of means-ends reasoning, in: {\em Common Sense, Reasoning, and Rationality} (ed. R. Elio),
Oxford Scholarship Online (2002)

\bibitem{rg95-1} Rao, M., Georgeff, M.: BDI-agents: from theory to practice, {\em  Proc. First Inter. Conf. Multiagent Syst.} (1995)

\bibitem{s55-1} Simon, H.: A behavioral model of rational choice, {\em Quarterly J. Econ.} {\bf 69}, 99--118 (1955)

\bibitem{scb88-1} Skinner, E., Chapman, M., Baltes, P.: Control, means-ends, and agency beliefs: a new conceptualization and its measurement during childhood,
{\em J. Personal. Social Psychol.} {\bf 54}, 117--133 (1988)

\bibitem{w00-1} Wooldridge, M.: {\em Reasoning About Rational Agents}, MIT Press (2000)
\end{thebibliography}

\begin{thebibliography}{99}

\bibitem{ae12-2} Aarts, H., Elliott, A.: {\em Goal-Directed Behavior}, Psychology Press (2012)

\bibitem{b88-2} Bandura, A.: Self-regulation of motivation and action through goal systems, in: {\em Cognitive Perspectives on Emotion and Motivation}
(eds. V. Hamilton, G.H. Bower, N.H. Frijda), pp.\ 37--41, Kluwer Academic Publishers (1988)

\bibitem{g97-2} Gutman, J.: Means--end chains as goal hierarchies, {\em Psychol. Market.} {\bf 14}, 545--560 (1997)

\bibitem{ksffcs02-2} Kruglanski, A., Shah, J., Fishbach, A., Friedman, R., Chun, W.Y., Sleeth-Keppler, D.: {\em A Theory of Goal Systems},  Chapter in: {\em
Advances in Experimental Social Psychology} (ed. M.P. Zanna) {\bf 34}, pp.\ 331--378, Academic Press (2002)

\bibitem{pba95-2} Pieters, R., Baumgartner, H., Alien, D.: A means-ends chain approach to consumer goal structures, {\em Intern. J. Res. Market.} {\bf 12},
227--244 (1995)

\bibitem{h97-2} Higgins, E.: Beyond pleasure and pain, {\em Amer. Psychol.} {\bf 52}, 1280--1300 (1997)

\bibitem{vw85-2} Vallacher, R., Wegner, D.: {\em A Theory of Action Identification},  Hillsdale (1985)

\bibitem{wm12-2} Wang, C., Mukhopadhyay, A.: The dynamics of goal revision: a cybernetic multiperiod test-operate-test-adjust-loop (TOTAL) model of
self-regulation, {\em J. Cons. Res.} {\bf 38}, 815--832 (2012)
\end{thebibliography}

\begin{thebibliography}{99}

\bibitem{a06-3} Andersen, E.: Schumpeter's general theory of social evolution: the early version, {\em Proc. Conf. on Neo-Schumpeterian Economics: An Agenda for
the 21st Century}, Trest, Czech Republic (2006)

\bibitem{a08-3} Andersen, E.: The essence of Schumpeter's evolutionary economics: a centennial appraisal of his first book,  {\em Proc. Intern.
Schumpeter Soc. Conf.} pp.\ 2--5, Rio de Janeiro, Brazil (2008)

\bibitem{a99-3} Argyris, C.: {\em On Organizational Learning}, Blackwell Publishers (1999)

\bibitem{b99-3} Bandura, A.: A social cognitive theory of personality, in: {\em Handbook of Personality} (eds. L. Pervin, O. John), pp.\ 154--196,
Guilford Publications (1999)

\bibitem{b93-3} Becker, G.: {\em Human Capital: A Theoretical and Empirical Analysis with Special Reference to Education}, 3rd edition,  University of Chicago
Press (1993)

\bibitem{bh87-3} Beckhard, R., Harris, R.: {\em Organizational Transitions: Managing Complex Change}, Addison-Wesley Publishers (1987)

\bibitem{b87-3} Bratman, M.: {\em Intentions, Plans, and Practical Reasons}, Harvard University Press (1987)

\bibitem{b91-3} Bridges, W.: {\em Managing Transitions}, Perseus Publishers (1991)

\bibitem{bm02-3} Bridges, W., Mitchell, S.: Leading transition: a new model for change, in: {\em On Leading Change} (eds. F. Hesselbein, R. Johnston),
Jossey-Bass Publishers (2002)

\bibitem{b80-3} Brisoux, J.: {\em Le Ph\'{e}nom\`{e}ne des Ensembles \'{e} Voqu\'{e}s: Une \'{E}tude Empirique des Dimensions Contenu et Taille},
Th\`{e}se de doctorat, Universit\'{e} Laval (1980)

\bibitem{bl81-3} Brisoux, J., Laroche, M.: Evoked set formation and composition: an empirical investigation under a routinized response behavior situation,
in: {\em Advances in Consumer Research} (ed. K.B. Monroe), Vol, 8, pp. 357--361, Ann Arbor, MI, Association for Consumer Research (1981)

\bibitem{bb85-3} Bullock, R., Batten, D.: It's just a phase we're going through, {\em Group Organiz. Stud.} {\bf 10}, 383--412 (1985)

\bibitem{c90-3} Carnall, C.: {\em Managing Change in Organizations}, Prentice Hall (1990)

\bibitem{cs98-3} Carver, C.S., Scheier, M.F.: {\em On the Self-Regulation of Behavior}, Cambridge University Press (1998)

\bibitem{cg04-3} Chen, Y, Gazzale, R.: When does learning in games generate convergence to Nash equilibria? The role of supermodularity
in an experimental setting, {\em Amer. Econ. Rev.} {\bf 5}, 1505--1535 (2004)

\bibitem{csos13-3} Cruz Neto, J.X., Oliveira, P.R., Soares Jr., P.A., Soubeyran, A.: Learning how to play Nash, potential games and alternating minimization
method for structured nonconvex problems on Riemannian manifolds, {\em J. Convex Anal.} {\bf 20}, 395--438 (2013)

\bibitem{rw06-3} De Ridder, D., de Wit, J.: Self-regulation in health behavior: concepts, theories, and central issues, in {\em Self-Regulation Health
Behavior} (eds. D. de Ridder, J. de Wit), pp.\ 1-23, John Wiley \& Sons (2006)

\bibitem{ec95-3} Eagly, A., Chaiken, S.: Attitude strength, attitude structure and resistance to change, in: {\em Attitude Strength} (eds. R. Petty, J. Kosnik),
pp.\ 413--432), Erlbaum Publishes (1995).

\bibitem{ebm95-3} Engel, J., Blackwell, R., Miniard, P.: {\em Consumer Behavior. Forth Worth}, Dryden Press (1995)

\bibitem{g90-3} Gollwitzer, P.: Action phases and mind-sets. Handbook of motivation and cognition, {\em Found. Social Behav.} {\bf 2}, 53--92 (1990)

\bibitem{ge77-3} Gould, S., Eldredge, N.: Punctuated equilibria: the tempo and mode of evolution reconsidered, {\em Paleobiolog.} {\bf 3}, 115--151 (1977)

\bibitem{g81-3} Gottfredson, L.: Circumscription and compromise: a developmental theory of occupational aspirations, {\em J. Counsel. Psychol.} {\bf 28},
545--579 (1981)

\bibitem{h07-3} Heckhausen, J.: The motivation-volition divide and its resolution in action-phase models of developmental regulation, {\em Res. Human
Develop.} {\bf 4}, 163--180 (2007)

\bibitem{h06-3} Hiatt, J.: ADKAR: A model for change in business, government and our community, {\em Prosci} (2006)

\bibitem{hs69-3} Howard, J.A., Sheth, J.N.: {\em The Theory of Buyer Behavior (Marketing)}, John Wiley \& Sons (1969)

\bibitem{jb84-3} Janz, N., Becker, M.: The health belief model: a decade later, {\em Health Educat. Behav.} {\bf 11}, 1--47 (1984)

\bibitem{k84-3} Kolb, D.: {\em Experiential Learning: Experience as the Source of Learning and Development}, Prentice-Hall (1984)

\bibitem{k95-3} Kotter, J.: Leading change: why transformation efforts fail?  {\em Harvard Business Review} (1995)

\bibitem{ksffcs02-3} Kruglanski, A., Shah, J., Fishbach, A., Friedman, R., Chun, W.Y., Sleeth-Keppler, D.: {\em A Theory of Goal Systems},  Chapter in: {\em
Advances in Experimental Social Psychology} (ed. M.P. Zanna) {\bf 34}, pp.\ 331--378, Academic Press (2002)

\bibitem{l47-3} Lewin, K.: Frontiers in group dynamics: concept, method and reality in social science, social equilibria and social change, {\em  Human Relat.}
{\bf 1}, 5--41 (1947)

\bibitem{l51-3} Lewin, K.: {\em Field Theory in Social Science}, Harper and Row Publishers (1951)

\bibitem{l59-3} Lindblom, C.: The science of ``muddling through", {\em Public Admin. Rev.}, 79--88 (1959)

\bibitem{lww58-3} Lippitt, R., Watson, J., Westley, B.: {\em Dynamics of Planned Change}, Harcourt and Brace Publushers (1958)

\bibitem{m91-3} March, J.: Exploration and exploitation in organizational learning, {\em Organiz. Sci.} {\bf 2}, 71--87 (1991)

\bibitem{m54-3} Maslow, A.: {\em Motivation and Personality}, Harper Publishers (1954)

\bibitem{mla85-3} McGuire, W., Lindzey, G., Aronson, E.: Attitudes and attitude change, {\em Handbook Soc. Psychol.} {\bf 2}, 233--346 (1985)

\bibitem{m38-3} Murray, H.: {\em Explorations in Personality}, Oxford University Press (1938)

\bibitem{nt80-3} Nadler, D., Tushman, M.: A model for diagnosing organizational behavior, {\em Organiz. Dynamics} {\bf 8}, 35--951 (1980)

\bibitem{nw82-3} Nelson, R., Winter, S.: {\em An Evolutionary Theory of Economic Change}, Harvard University Press (1982)

\bibitem{pd92-3} Prochaska, J., DiClemente, C.: Stages of change in the modification of problem behaviors, {\em Prog. Behav. Modif.} {\bf 28}, 183--218 (1992)

\bibitem{pd94-3} Prochaska, J., DiClemente, C.: {\em Changing for Good: the Revolutionary Program that Explains the Six Stages of Change and Teaches You
How to Free Yourself from Bad Habits}, Morrow Publishers (1994)

\bibitem{s99-3} Schein, E.: {\em The Corporate Culture Survival Guide: Sense and Nonsense about Culture Change}, Jossey-Bass Publishers (1999)

\bibitem{s11-3} Schoenfeld, A.: {\em How We Think: A Theory of Human Decision-Making with the Focus on Teaching}, Routledge Publishers (2011)

\bibitem{s02-3} Schumpeter, J.: New translations from theorie der wirtschaftlichen entwicklung, {\em  Amer. J. Econ. Sociol.} {\bf 61}, 405--437 (2002)

\bibitem{s05-3} Schumpeter, J.: On the nature of economic crises, in: {\em Business Cycle Theory: Selected Texts 1860--1939}  (ed. M. Boianovsky),
Pickering \& Chatto, Vol.\ 5, pp.\ 5--50 (2005)

\bibitem{s90-3} Senge, P.: {\em The Fifth Discipline: The Art and Practice of the Learning Organization}, Doubleday Publishers (2005)

\bibitem{si44-3} Simpson, G.G.: {\em Tempo and Mode in Evolution}, Columbia University Press (1944)

\bibitem{s94-3} Snyder, C.: {\em The Psychology of Hope: You Can Get There from Here}, New York Free Press (1994)

\bibitem{sbcg13} Soares, C., Batista, L., Campelo, F., Frederico Gadelha Guimaraes, F.G.:  Computation of mixed strategy nondominated Nash equilibria in
game theory,  {\em Proc. of the 1st BRICS Countries Congress}, Brazil (2013)

\bibitem{wn07-3} Wood, W., Neal, D.: A new look at habits and the habit--goal interface, {\em Psychol. Rev.} {\bf 114}, 843--863 (2007)

\bibitem{w02-3} Wooldridge, M.: {\em An Introduction to Multi-Agent Systems}, John Wiley \& Sons (2002)
\end{thebibliography}

\begin{thebibliography}{99}

\bibitem{a90-4} Abel, A.: Asset prices under habit formation and catching up with the joneses, {\em  Amer. Econ. Rev.} {\bf 80}, 38--42 (1990)

\bibitem{afgh11-4} Abeler, J., Falk, A.,  G\"{o}tte, L., Huffman, D.: Reference points and effort provision, Discussion Paper No.\ 358, GESY (2011)

\bibitem{akb09-4} Abdeljawad, T., Karapinar, E., Ben-El-Mechaiekh, H.: Quasicone metric spaces and generalizations of Caristi-Kirk's theorem, {\em  J. Fixed
Point Theory Appl.} {\bf 13}, 574--387 (2009)

\bibitem{a04-4} Appadurai, A.: The capacity to aspire: culture and the terms of recognition, in: {\em Culture and Public Action} (eds. V. Rao, M. Walton),
 pp.\ 59--84, Stanford University Press (2004)

\bibitem{bc13-4} Bae, J., Cho, S.: Fixed points and variational principle with applications to equilibrium problems on cone metric spaces, {\em  J. Korean Math.
Soc.} {\bf 50}, 95--109 (2013)

\bibitem{bs81-4} Bandura, A., Schunk, D.: Cultivating competence, self efficacy, and intrinsic interest through proximal self motivation, {\em  J. Personal.
Soc. Psychol.} {\bf 41}, 586--598 (1981)

\bibitem{bm07-4} Bao, T.Q., Mordukhovich, B.S.: Variational principles for set-valued mappings with applications to multiobjective optimization, {\em
Control Cyber.} {\bf 36}, 531--562 (2007)

\bibitem{bm10-4} Bao, T.Q., Mordukhovich, B.S.: Relative Pareto minimizers for multiobjective problems: existence and optimality conditions, {\em Math. Program.}
{\bf 122}, 301--347 (2010)

\bibitem{bm13-4} Bao, T.Q., Mordukhovich, B.S.: Necessary nondomination conditions for set and vector optimization with variable structures,
{\em J. Optim. Thepry Appl.}, DOI 10.1007/s10957-013-0332-6 (2013)

\bibitem{b94-4} Bargh, J.: The four horsemen of automaticity: awareness, intention, efficiency, and control in social cognition,  {\em Handbook of Social
Cognition}, Vol.\ 1, pp.\ 1-40, Erlbaum Publishers (1994)

\bibitem{bm88-4} Becker, G., Murphy, K.: A theory of rational addiction, {\em J. Polit. Econ.} {\bf 96}, 675--700 (1988)

\bibitem{bb13-4} Bello Cruz, J., Bouza Allende, G.: A steepest descent-like method for variable order vector optimization problems, {\em  J. Optim.
Theory Appl.}, to appear (2013)

\bibitem{b13-4} Bhatia, S.: Associations and the accumulation of preference, {\em Psychol. Rev.} {\bf 120}, 522--543 (2013)

\bibitem{bg12-4} Bhatia, S., Golman, R.: Attention and reference dependence, Department of Social \& Decision Sciences, Carnegie Mellon University (2012)

\bibitem{b87-4} Bratman, M.: {\em Intentions, Plans, and Practical Reasons}, Harvard University Press (1987)

\bibitem{c00-4} Carroll, C.: Solving consumption models with multiplicative habits, {\em Econ. Letts.} {\bf 68}, 67--77 (2000)

\bibitem{cy02-4} Chen, G., Yang, X.: Characterizations of variable domination structures via nonlinear scalarization,  {\em J. Optim. Theory Appl.}
{\bf 112}, 97--110 (2002)

\bibitem{cb12-4} Cho, S., Bae, J.: Variational principles on cone metric spaces, {\em Int. J. Pure Appl. Math.} {\bf 77}, 709--718 (2012)

\bibitem{cjcnw04-4} Cooke, A., Janiszewski, C., Cunha, M., Nasco, S., De Wilde, E.: Stimulus context and the formation of consumer ideals, {\em J. Consum. Res.}
{\bf 31}, 112--124 (2004)

\bibitem{e05-4} Easterlin, R.: A puzzle for adaptive theory, {\em J. Econ. Behav. Organiz.} {\bf 56}, 513--521 (2005)

\bibitem{e11-4} Eichfelder, G.: Optimal elements in vector optimization with a variable ordering structure, {\em  J. Optim. Theory Appl.} {\bf 151},
217--240 (2011)

\bibitem{e13-4} Eichfelder, G.: Ordering structures in vector optimization and applications in medical engineering, Technical University Ilmenau, Institute
of Mathematics, Preprint (2013)

\bibitem{eh13-4} G. Eichfelder and T.X.D. Ha, Optimality conditions for vector optimization problems with variable ordering structures,
{\em Optimization}, to appear (2013).

\bibitem{gl11-4} Golman, R., Loewenstein, G.: Explaining non convex preferences with aspirational and status quo reference dependence, Mimeo, Carnegie Mellon
University (2011)

\bibitem{grtz03-4} G\"{o}pfert, A., Riahi, H., Tammer, C., Z\u{a}linescu, C.: {\em Variational Methods in Partially Ordered Spaces},
Springer (2003)

\bibitem{ge77-4} Gould, S., Eldredge, N.: Punctuated equilibria: the tempo and mode of evolution reconsidered, {\em Paleobiol.} {\bf 3}, 115--151 (1977)

\bibitem{gh00-4} Granas, A., Horvath, C.: On the order theoretical Cantor theorem, {\em Taiwanese J.of Math.} {\bf 4}, 203--213 (2000)

\bibitem{gjn08-4} Gutierrez, C., Jimenez, B., Novo, V.: A set-valued Ekeland's variational principle in vector optimization, {\em  SIAM J. Control Optim.} {
\bf 47}, 883--903 (2008)

\bibitem{hw90-4} Hauser, J., Wernerfelt, B.: An evaluation cost model of consideration sets, {\em J. Consum. Res.} {\bf 16}, 393--408 (1990)

\bibitem{hlw99-4} Heath, C., Larrick, R., Wu, G.: Goals as reference points, {\em Cognit. Psychol.} {\bf 38}, 79--109 (1999)

\bibitem{h64-4} Helson, H.: Current trends and issues in adaptation-level theory, {\em Amer. Psycholog.} {\bf 19}, 26--38 (1964)

\bibitem{kt79-4} Kahneman, D, Tversky, A.: Prospect theory: an analysis of decision under risk, {\em Economet.} {\bf 47}, 263--291 (1979)

\bibitem{kq11-4} Khanh, P.Q., Quy, D.N.: On generalized Ekeland's variational principle and equivalent formulations for set-valued mappings, {\em J. Global
Optim.} {\bf 49}, 381--396 (2011)

\bibitem{ko87-4} Klein, N., Oglethorpe, J.: Cognitive reference points in consumer decision making, in: {\em Advances in Consumer Research} (eds.
M. Wallendorf, P. Anderson), Vol. 14, pp.\ 183--187, Association for Consumer Research (1987)

\bibitem{l59-4} Lindblom, C.: The science of ``muddling through", {\em Public Admin. Rev.}, 79--88 (1959)

\bibitem{ls13-4} Luc, D.T., Soubeyran, A.: Variable preference relations: existence of maximal elements, {\em J. Math.  Econ.}, to appear (2013)

\bibitem{lm11-4} Liu, C. G., Ng, K.F.: Ekeland's variational principle for set-valued functions, {\em SIAM J. Optim.} {\bf 21} 41--56 (2011)

\bibitem{q12-4} Qiu, J.H.: On Ha's version of set-valued Ekeland's variational principle, {\em  Acta Math. Sinica (English Series)} {\bf 28}, 717--726 (2012)

\bibitem{rsu06-4} Ravn, F., Schmitt-Groh\'{e}, S., Uribe, M.: Deep habits, {\em Rev. Econ. Stud.} {\bf 73}, 1--24 (2006)

\bibitem{r75-4} Rosch, E.: Natural categories, {\em Cognit. Psychol.} {\bf 4}, 328--350 (1975)

\bibitem{s88-4} Selten, R.: Aspiration adaptation theory, {\em J. Math. Psychol.} {\bf 42}, 210--214 (1988)

\bibitem{sn13-4} Shaddad, F., Noorani, M.: Fixed point results in quasi-cone metric spaces, {\em Abstract Appl. Anal.} (2013)

\bibitem{s55-4} Simon, H.: A behavioral model of rational choice, {\em Quarterly J. Econ.} {\bf 69}, 99--118 (1955)

\bibitem{s94-4} Snyder, C.: {\em The Psychology of Hope: You Can Get There from Here}, New York Free Press (1994)

\bibitem{t99-4} Thaler, R.: Mental accounting matters, {\em J. Behav. Dec. Making} {\bf 12}, 183--206 (1999)

\bibitem{tk91-4} Tversky, A., Kahneman, D.: Loss aversion in riskless choice: a reference-dependent model, {\em Quarterly J. Econ.} {\bf 106},
1039--1061 (1991)

\bibitem{wj06-4} Weber, E., Johnson, E.: Constructing preferences from memories, in {\em The Construction of Preference} (eds. S. Lichtenstein,
P. Slovic), pp.\ 397--410, Cambridge University Press (2006)

\bibitem{w02-4} Wendner, R.: Habits: multiplicative or subtractive? Depart. Econ., Graz University (2002)

\bibitem {wn07-4} Wood, W., Neal, D.: A new look at habits and the habit--goal interface, {\em Psychol. Rev.} {\bf 114}, 843--863 (2007)

\bibitem{w00-4} Wooldridge, M.: Reasoning about rational agents,  MIT Press (2000)

\bibitem{wo02-4} Wooldridge, M.: {\em An Introduction to Multi-Agent Systems}, John Wiley \& Sons (2002)

\end{thebibliography}
\end{document}